\documentclass{lms}
\usepackage{amssymb}
\usepackage{amsmath}
\usepackage{verbatim}
\usepackage{color}
\usepackage{url}

\def\un{{\mathrm{1~\hspace{-1.4ex}l}}}

\newcommand{\BigO}{\mathcal{O}}
\newcommand{\jvRe}{\operatorname{Re}}
\newcommand{\jvIm}{\operatorname{Im}}

\newcommand{\eps}{\varepsilon}
\newcommand{\frakp}{\mathfrak{p}}
\newcommand{\frakq}{\mathfrak{q}}
\newcommand{\frakr}{\mathfrak{r}}
\newcommand{\frakU}{\mathfrak{U}}
\newcommand{\opnm}{\operatorname}

\newtheorem{theorem}{Theorem}[section]
\newtheorem{proposition}[theorem]{Proposition}

\newtheorem{lemma}[theorem]{Lemma}

\newnumbered{remark}{Remark}
\newnumbered{example}{Example}
\newnumbered{definition}{Definition}

\numberwithin{equation}{section}

\author{Joe Viola}
\title[Resolvent Estimates for NSA Operators]{Resolvent Estimates for Non-Selfadjoint Operators with Double Characteristics\footnote[1]{J. London Math. Soc., to appear.  Published version may differ.}}

\classno{47A10 (primary), 35P15 (secondary)}

\begin{document}

\maketitle

\begin{abstract}

We study resolvent estimates for non-selfadjoint semiclassical pseudodifferential operators with double characteristics.  Assuming that the
quadratic approximation along the double characteristics is elliptic, we obtain polynomial upper bounds on the resolvent in a suitable region
inside the pseudospectrum.

\end{abstract}

\section{Introduction and Statement of Results}

\subsection{Related Prior Results Bounding Resolvents}

The spectrum of a differential operator on $L^2(\Bbb{R}^d)$ given by quantizing an elliptic quadratic form $q: \Bbb{R}^d_x\times \Bbb{R}^d_\xi
\rightarrow \Bbb{C}$ has been known since the work of J. Sj\"{o}strand in 1974 \cite{Sj74}.  Here, a quadratic form $q$ is said to be elliptic when
\begin{equation}\label{eDefEllipticQuad}
q(x,\xi) = 0 \Leftrightarrow (x,\xi) = 0, \qquad q(\Bbb{R}^{2d}) \neq \Bbb{C}.
\end{equation}
It is proved in~\cite{Sj74} that the range of $q$ is either a closed cone lying properly within a half-plane of $\Bbb{C}$ or all of $\Bbb{C}.$ The latter is a degenerate case which can only happen, assuming $q^{-1}(\{0\}) = \{(0,0)\},$ when $d = 1,$ and it will not be considered here.
The quantization we will use throughout is the $h$-dependent Weyl quantization, and for quadratic forms
$$
	q(x,\xi) = \sum_{|\alpha| + |\beta| = 2}q_{\alpha \beta}x^\alpha\xi^\beta,
$$
this takes the form
\begin{equation}\label{eExplicitQuadraticWeyl}
	q^w(x,hD_x) = \sum_{|\alpha| + |\beta| = 2} q_{\alpha \beta}\left(\frac{x^\alpha (hD_x)^\beta + (hD_x)^\beta x^\alpha}{2}\right),\quad 0< h \leq 1.
\end{equation}
Derivatives are written using the standard notation, $D_{x_j} = \frac{1}{i}\partial_{x_j}.$

This paper is concerned with resolvent bounds for pseudodifferential operators.  For completeness, we recall the basic definitions of spectrum and
resolvent. The spectrum of a closed densely defined operator $A$ on a Hilbert space $\mathcal{H}$ is defined via its complement,
$$\operatorname{Spec}(A)^c = \{z \in \Bbb{C} \::\: (A - z):\operatorname{Dom}(A)\rightarrow \mathcal{H}~\mathrm{bijective}\},$$ and the resolvent
for $z \notin \operatorname{Spec}(A)$ is simply $(A-z)^{-1}.$  The resolvent norm will refer to the standard operator norm
\begin{eqnarray*}
	||(A-z)^{-1}||_{\mathcal{H}\rightarrow\mathcal{H}}
		&=& \sup_{x \in \mathcal{H}, ||x|| = 1} ||(A-z)^{-1}x||
		\\ &=& \left(\inf_{x \in \operatorname{Dom}(A), ||x|| = 1} ||(A-z)x||\right)^{-1}.
\end{eqnarray*}

In \cite{Sj74}, we have a complete description of the spectrum of differential operators given by quantizing quadratic forms which are elliptic in the
sense of (\ref{eDefEllipticQuad}). When recalling this description, let us assume, in order to fix the ideas, that $\jvRe q$ is positive definite. With
$$
	\sigma\left((x,\xi),(y,\eta)\right) = \langle \xi, y \rangle - \langle \eta, x \rangle = \sum_{j=1}^n (\xi_j y_j - \eta_j x_j)
$$
the standard symplectic product on $\Bbb{C}^{2d},$ we recall that the Hamilton map $F$ of a quadratic form $q(x,\xi)$ is the unique $2d \times 2d$ matrix with entries in $\Bbb{C}$ for which
$$
	q(x,\xi) = \sigma\left((x,\xi),F(x,\xi)\right),\quad \sigma\left((x,\xi), F(y,\eta)\right) = -\sigma\left(F(x,\xi),(y,\eta)\right)
$$
for every $(x,\xi), (y,\eta) \in \Bbb{R}^{2d}$.

Having established this notation, we recall that, according to Theorem 3.5 of~\cite{Sj74}, the spectrum of $q^w(x,hD_x)$ is composed
of eigenvalues of finite multiplicity, determined by the eigenvalues of the Hamilton map $F$ of $q$:
\begin{equation}\label{eSjExplicitSpec}
	\operatorname{Spec}(q^w(x,hD_x)) = \left\{ \mathop{\sum_{\lambda \in \operatorname{Spec}F,}}_{\jvIm\lambda > 0}\frac{h}{i}(r_\lambda + 2k_\lambda)\lambda\::\: k_\lambda \in \Bbb{N}\cup \{0\} \right\}.
\end{equation}
Here, $r_\lambda$ is the dimension of the generalized eigenspace of $F$ which corresponds to $\lambda.$

For $z \in \Bbb{C}$ with $|z| \sim h,$ the growth of the operator norm of the resolvent $(q^w(x,hD_x)-z)^{-1}$, as $h\rightarrow 0^+$,
is at most a constant times $h^{-1},$ whether $q^w(x,hD_x)$ is a normal operator or not.  Indeed, there is a change of variables which shows that
$q^w(x,hD_x)$ is unitarily equivalent to $hq^w(x,D_x),$ so if $q^w$ is self-adjoint (that is, $q(x,\xi)$ is real-valued) and
$K \subset \Bbb{C}$ is a fixed compact set containing no eigenvalues of $q^w(x,D_x),$ we have
\begin{eqnarray*}
	||(q^w(x,hD_x) - hz)^{-1}||_{L^2(\Bbb{R}^d)\rightarrow L^2(\Bbb{R}^d)} &=& \frac{1}{\operatorname{dist}(hz, \operatorname{Spec}q^w(x,hD_x))}
	\\ &=& h^{-1}\frac{1}{\operatorname{dist}(z, \operatorname{Spec}q^w(x,D_x))},
	\\ &=& \BigO_K(h^{-1}),
\end{eqnarray*}
for all $z \in K.$  On the other hand, if $q^w(x,hD_x)$ is not normal, one can simply take the supremum of $||(q^w(x,D_x) - z)^{-1}||$ over $z \in K$.
Then the scaling change of variables mentioned above gives the same qualitative result,
\begin{equation}\label{IntroScaledQuadratic}
	||(q^w(x,hD_x) - hz)^{-1}||_{L^2(\Bbb{R}^d)\rightarrow L^2(\Bbb{R}^d)} = \BigO_K (h^{-1}).
\end{equation}

When $|z|\gg h,$ however, recent work involving non-selfadjoint harmonic oscillators (see E. B. Davies and A. B. J. Kuijlaars \cite{Davies04}, L. Boulton
\cite{Boulton02}) showed that rapid growth of the resolvent is possible, even if $z$ maintains some distance from
$\operatorname{Spec}(q^w(x,hD_x))$ while remaining in $q(\Bbb{R}^{2d}).$

In the specific context of elliptic quadratic differential operators, K. Pravda-Starov \cite{PSDuke} proved that the resolvent of $q^w(x,hD_x)$ grows superpolynomially for every point in the interior of the range, assuming that $q^w(x,hD_x)$ fails to be normal.  This assumption holds exactly when the Poisson bracket of real and imaginary parts of the symbol is nonzero:
$$
	\{\jvRe q, \jvIm q\} = \frac{\partial \jvRe q}{\partial \xi}\cdot \frac{\partial \jvIm q}{\partial x} - \frac{\partial \jvIm q}{\partial \xi}\cdot \frac{\partial \jvRe q}{\partial x} \not\equiv 0
$$
on $\Bbb{R}^{2d}.$  Under this hypothesis, for every $z \in (q(\Bbb{R}^{2d}))^\circ$ and $N \in \Bbb{N},$ there exist $h_0 > 0$ and a family $\{u_h\}_{h \in (0,h_0]}$ in $\mathcal{S}(\Bbb{R}^d)$ where $||u_h||_{L^2} = 1$ while, for $h \in (0, h_0],$
$$
	||(q^w(x,hD_x)-z)u_h||_{L^2} = \BigO(h^N).
$$
It follows that
$$
	||(q^w(x,hD_x) - z)^{-1}|| \geq \frac{h^{-N}}{C_N} \quad \forall N \in \Bbb{N},
$$
as $h \rightarrow 0^+,$ unless $z \in \operatorname{Spec}(q^w(x,hD_x)).$ In fact, because the symbol is analytic, the lower bound for the resolvent grows exponentially in $1/h$,
which follows from the work by N. Dencker, J. Sj\"{o}strand, and M. Zworski \cite{DeSjZw04}.

Studies of bounds on the operator norm of the resolvent of a non-selfadjoint operator are equivalent to the study of the
pseudospectrum of that operator. The $\varepsilon$-pseudospectrum of a closed densely defined operator $A$, acting on some Hilbert space, is defined as
$$
	\operatorname{Spec}_\varepsilon(A) = \{z \in \Bbb{C} \::\: ||(A-z)^{-1}|| > \varepsilon^{-1}\}.
$$
We refer the reader to the recent monograph \cite{TrefEmb05} by L. N. Trefethen and M. Embree for an extensive discussion of
pseudospectra with many applications.  In the context of $h$-pseudodifferential operators, the term "semiclassical pseudospectrum" was used in
\cite{DeSjZw04} to make precise the relationship between a subset of $\Bbb{C}$ where the resolvent grows rapidly as $h\rightarrow 0^+$,
and a subset of the range of the symbol.  Other semiclassical versions of the pseudospectrum appear in \cite{PSDuke}.
In the latter paper, as in \cite{H-PS}, the definition of semiclassical pseudospectrum uses a fixed spectral parameter $z\in \Bbb{C}$, while in the present work the spectral parameter is is allowed to vary only in $h$-dependent regions.  We therefore make no precise analogy with these definitions for semiclassical pseudospectrum.  The motivation, however, comes from these earlier works.

A natural extension from the quadratic case to pseudodifferential operators
$$
	p^w(x,hD_x)u(x) =\frac{1}{(2\pi h)^d} \iint_{\Bbb{R}^{2d}} e^{\frac{i}{h}(x-y)\cdot \xi}p(\frac{x+y}{2},\xi)u(y)\,dy\,d\xi
$$
is to consider symbols $p \in C_b^\infty(\Bbb{R}^{2d},\Bbb{C})$, bounded with all derivatives, with finitely points minimizing $\jvRe p$ and elliptic elsewhere.  At these points we assume the symbol is approximated to second order by a quadratic form with elliptic real part.  Without loss of generality it may be assumed that the
minimum of $\jvRe p$ is zero and occurs at precisely one doubly characteristic point, the origin, where we have a quadratic approximation with elliptic real part.  This is because a different minimum realized at a different $(x_0,\xi_0)\in\Bbb{R}^{2d}$ may be reduced to this case by elementary operations, and if finitely many points realize the minimum of $\jvRe p$, each may be considered separately and the results may be combined using localization via smooth cutoff functions.

For many such operators, even when non-selfadjoint, it is known that a resolvent bound of type (\ref{IntroScaledQuadratic}) holds.  Writing $q$ for the quadratic part of $p$ at the doubly characteristic part, if
$K \subseteq \Bbb{C}$ is compact and $K \cap (\opnm{Spec}q^w(x,D_x)) = \varnothing$, we have,
for all sufficiently small $h,$ the resolvent bound
$$
	||(p^w(x,hD_x) - hz)^{-1}||_{L^2(\Bbb{R}^d)\rightarrow L^2(\Bbb{R}^d)} = \BigO_K(h^{-1}).
$$
In the case when the quadratic approximation is elliptic,
this has been essentially well-known since \cite{Sj74}. In cases more general than the elliptic case, where the quadratic part is allowed to be
zero in certain directions but is subject to some additional nondegeneracy
assumptions, the same result has been
successfully shown in such recent papers as \cite{HeSjSt05} and \cite{H-PS}.

To the author's knowledge, results pushing the set from which the spectral parameter $z$ is taken beyond a region of size $\BigO(h)$ are heretofore unknown. (For related
results in the non-elliptic case, see \cite{HeSjSt05}.) Because the resolvent of an operator with elliptic quadratic symbol is known to grow rapidly,
some loss in resolvent norm beyond $\BigO(h^{-1})$ is necessary. However, in order to use standard methods joining together estimates for the quadratic approximation with
estimates from ellipticity for $(x,\xi)$ away from zero, it is necessary to keep the additional loss to only a small power of $h.$  This paper explores
sufficient assumptions on the values of $|z| \gg h$ and the symbol $p$ for which we have resolvent estimates of the form $$||(p^w(x,hD_x) -
z)^{-1}||_{L^2(\Bbb{R}^d)\rightarrow L^2(\Bbb{R}^d)} = \BigO(h^{-1-\gamma}),$$ where $\gamma > 0$ is taken sufficiently small but fixed,
as $h \rightarrow 0^+.$

\subsection{Hypotheses on Symbol}\label{SymbolHypotheses}

The class of symbols considered here will be bounded with all
derivatives: 
$$
	p(x,\xi) \in S(1) := \{a\in C^\infty(\Bbb{R}^{2d},\Bbb{C})\::\: |\partial^\alpha a| = \BigO_\alpha(1) ~\forall \alpha \in \Bbb{N}^{2d}\}.
$$
Furthermore, we will assume that the real part of our symbol will be
nonnegative, $\jvRe p \geq 0$; zero precisely at the origin,
$(\jvRe p)^{-1}(\{0\}) = \{(0,0)\}$; and satisfying
$$
	\liminf_{|(x,\xi)|\rightarrow \infty} \jvRe p > 0.
$$
As a consequence, we are assuming that $\jvRe p(x,\xi) > 1/C_r$ on any set $\{|(x,\xi)| > r > 0\}$.

\begin{remark}\label{RemFredholmFamily}

General operator-theoretic considerations (outlined in Section \ref{SectionGeneralUpperBounds}) show that $$p^w(x,hD_x) - z= \BigO(1):
L^2(\Bbb{R}^d) \rightarrow L^2(\Bbb{R}^d)$$ is an analytic family of Fredholm operators of index 0 for $z \in \operatorname{neigh}(0, \Bbb{C})$,
for $h$ sufficiently small.  In particular, injectivity and invertibility are equivalent for such $z.$

\end{remark}

At the origin, we assume that our symbol has a doubly characteristic point, $$p(0,0) = p'(0,0) = 0.$$  Consequently, $$p(x,\xi) = q(x,\xi) +
\BigO(|(x,\xi)|^3),$$ for $q(x,\xi)$ a quadratic form with $\jvRe q \geq 0.$ The final assumption we make is that
$\jvRe q(x,\xi)$ is an elliptic quadratic form in the sense of (\ref{eDefEllipticQuad}).  As a consequence,
$$\frac{1}{C}|(x,\xi)|^2 \leq \jvRe q(x,\xi) \leq C|(x,\xi)|^2$$ for
some constant $C > 0.$

\begin{example}

The hypotheses above are satisfied by a Schr\"{o}dinger operator with complex potential with some assumptions on the potential.  Specifically, let us
consider
$$p(x,\xi) = \chi(\xi)\xi^2 + (1-\chi(\xi)) + V(x) + iW(x),$$ where $V, W\in S(1)$ as functions from $\Bbb{R}^d$ to $\Bbb{R}$ and $\chi(\xi)\in
C_0^\infty(\Bbb{R}^d)$ is a standard cutoff function taking values in $[0,1]$ and equal to 1 in a neighborhood of $0 \in \Bbb{R}^d$. Furthermore assume $V(x) \geq 0$ and that both $V^{-1}(\{0\}) = \{0\}$
and $V(x) > 1/C$ for $|x| > C>0.$ Near zero, $$V(0) = W(0) = V'(0) = W'(0) = 0,$$ and $V''(0) > 0$ in the sense of quadratic forms.  Then all the
hypotheses on the symbol are met, and the theorem below applies.

\end{example}

The multiplication by a cutoff function in the $\xi$-variable to make the symbol in $S(1)$ does not reduce the applicability of the theorem.
The end of Section 1 in the paper \cite{H-PS} by M. Hitrik and K. Pravda-Starov shows that it is possible to extend the resolvent estimate from
symbols in $S(1)$ to those in $S(m)$, for $m \geq 1$ a general order function, as long as a natural assumption of ellipticity at infinity is satisfied.

\subsection{Main Theorem and Structure of Paper}

For $p \in S(1)$, a symbol satisfying the hypotheses given in Subsection \ref{SymbolHypotheses}, we have a simple resolvent estimate for
$z \in \Bbb{C}$ for which $h \ll |z| \leq hf(h),$ for $f(h)$ defined below. In order to limit the price paid for the growth factor
$f(h)$ to a resolvent growing no faster than $\BigO(h^{-1-\gamma})$ with $\gamma > 0$ small, we define
\begin{equation} \label{Definingf}
	f(h) = \frac{1}{C_\gamma}\left(\frac{\log\frac{1}{h}}{\log\log\frac{1}{h}}\right)^{1/d}.
\end{equation}
We will take $C_\gamma = C'/\gamma^{1/d}$ for $C'>0$ a geometric constant depending solely on the quadratic part of $p$ at the doubly characteristic point, but we are only free to do so when
$\gamma$ remains sufficiently small. With this choice of $f(h),$ we establish the following theorem.

\begin{theorem}\label{MainTheorem}

Let $p \in S(1)$ be a symbol satisfying the hypotheses in section \ref{SymbolHypotheses}.  For any $\gamma \in (0, 1/8),$ define $f(h)$ as in
(\ref{Definingf}). Then there exists $h_0 > 0$ such that, for $0 < h \leq h_0$ and $|z|\leq hf(h)$ with $$\operatorname{dist}(z, \operatorname{Spec}
(q^w(x,hD_x))) \geq hf(h)^{(1-d)/2},$$ the resolvent $(p^w(x,hD_x)-z)^{-1}: L^2(\Bbb{R}^d)\rightarrow L^2(\Bbb{R}^d)$ exists and obeys the bound
\begin{equation}\label{EqnMainResolventBound}||(p^w(x,hD_x) - z)^{-1}|| \leq \BigO(h^{-1-\gamma}).\end{equation}

\end{theorem}

A natural question is whether we could extend our set of allowed $z$ to a larger set of the form $\{|z| \leq h^\delta\}$, for $\delta \in (0,1).$ A
demonstration of useful scaling techniques shows that this cannot be the case.

Take any elliptic quadratic form $q(x,\xi)$ with $q(\Bbb{R}^{2d}) \neq \Bbb{C}$ and obeying $$\{\jvRe q, \jvIm q\} \not\equiv 0.$$  The scaling outlined in Appendix
\ref{ChangeVarsSubsection} gives that $$q^w(x,hD_x) = h^\delta q^w(\tilde{x}, \frac{h}{h^\delta}D_{\tilde{x}})$$ for $\tilde{x} = h^{-\delta/2}x.$
Let us write $\tilde{h} = h^{1-\delta}.$ Then by K. Pravda-Starov \cite{PSDuke}, for every $z \in (q(\Bbb{R}^{2d}))^\circ$ and for every $N \in \Bbb{N},$ there
exists $\tilde{h}_0 > 0$ and a family $\{u_{\tilde{h}}(\tilde{x})\}_{\tilde{h} \in (0, \tilde{h}_0]}$ in $\mathcal{S}(\Bbb{R}^d)$ with
$||u_{\tilde{h}}||_{L^2} = 1$ and $$||(q^w(\tilde{x},\tilde{h}D_{\tilde{x}}) - z)u_{\tilde{h}}||_{L^2} = \BigO(\tilde{h}^N).$$

Performing a unitary change of variables to change $\tilde{x}$ to $x$ and $\tilde{h}$ to $h$ gives $$||(q^w(x,hD_x) - h^\delta
z)U_{h^\delta}u_{\tilde{h}}|| = \BigO(h^{\delta + (1-\delta)N}).$$  Here $U_{h^\delta}u_{\tilde{h}}(\tilde{x}) = h^{\delta d/4}u_{\tilde{h}}(x)$ is
the unitary change of variables as in Appendix \ref{ChangeVarsSubsection}.

This demonstrates that if $|z| \sim h^\delta$ for $\delta \in (0, 1)$ and $z \in (q(\Bbb{R}^{2d}))^\circ,$ then the resolvent norm of
$q^w(x,hD_x)$ at $z$ grows superpolynomially quickly, and polynomial resolvent bounds are impossible.

The plan of this paper is as follows.  In Section \ref{SectionGeneralUpperBounds}, we derive essentially well-known general upper bounds on
resolvents of $h$-pseudodifferential operators, which require only weak hypotheses.  In Section \ref{SectionQuadraticCase}, those general upper
bounds are applied to a rescaled version of the approximating quadratic operator.  Section \ref{SectionLocalizing} introduces a cutoff function on the FBI transform side
localizing to a small $h$-dependent neighborhood of the origin on which the full symbol and its quadratic part are close together.  The region
outside that neighborhood is dealt with in Section \ref{SectionExterior}, and these results are combined to prove the main theorem in Section
\ref{SectionGluingProof}. Finally, Appendix \ref{ChangeVariablesAppendix} outlines the FBI-Bargmann transform and useful scaling relations, and
Appendix \ref{WeightsAppendix} discusses the equivalence of the graph norms and domains of elliptic quadratic operators, both on the real and on
the FBI-Bargmann sides.

\medskip

\begin{acknowledgements} The author would like to thank M. Hitrik for his very much appreciated help during the creation of this paper.
Furthermore, the author would like to thank N. Trefethen for suggestions improving the exposition in the introduction.  Finally, the author
would like to extend his gratitude to the referee for many detailed suggestions and corrections, which he feels have greatly improved the paper.
\end{acknowledgements}

\section{General Upper Bounds on Resolvents}\label{SectionGeneralUpperBounds}

For completeness, we include a proof, using well-known methods (see \cite{MarkusPencils}, \cite{Sj01}, \cite{DeSjZw04}) of an exponential dimension-dependent
upper bound for the resolvent of pseudodifferential operators elliptic at infinity whose range avoids a point.  The proof is of a general
operator-theoretic nature and disregards any special features of quadratic differential operators, which have been described in some detail in
\cite{Sj74}, and so improvements taking advantage of these features may be anticipated. However, those improvements could not give better-than-exponential upper bounds
thanks to the quasimodes found in the works \cite{PSDuke} and \cite{DeSjZw04}, as mentioned in the introduction. Future improvements to this proposition also could be inserted with minimal
difficulty in the rest of the proof, which could lead to an improvement in the growth factor $f(h)$ described in (\ref{Definingf}).

In order to apply this result to the unbounded quadratic part $q(x,\xi) \notin S(1),$ we consider symbols in standard symbol classes
$$S(m) := \{a \in C^\infty(\Bbb{R}^{2d})\::\:
|\partial_x^\alpha \partial_\xi^\beta a(x,\xi)| \leq \BigO_{\alpha \beta}m(x,\xi)\}$$ for any order function $m\geq 1$ fulfilling
$$m(x,\xi) \leq C_0 (1 + (x-y)^2 + (\xi - \eta)^2)^{N_0}m(y,\eta)$$ for some fixed $C_0, N_0 > 0.$
We shall assume here, as we may, that $m\in S(m)$, and introduce, for $h>0$ small enough, the natural Sobolev space associated with the order
function $m$,
$$
H(m) = \left(m^w(x,hD_x)\right)^{-1}(L^2(\Bbb{R}^d)).
$$

\begin{proposition}\label{ResolventUpperBounds}
Let $m \geq 1$ be an order function as above.  Consider $p \in S(m)$ for which $\jvRe p(x,\xi) \geq 0$ for all $(x,\xi)\in\Bbb{R}^{2d}.$
Assume in addition that $\jvRe p$ is elliptic at infinity in the sense that
\begin{equation}\label{eEllinfty}
	\jvRe p(x,\xi) \geq \frac{m(x,\xi)}{C}, \quad\quad \forall |(x,\xi)| \geq \tilde{C}
\end{equation}
for some $C, \tilde{C} > 0$ sufficiently large.  Let $\rho, \rho'$ be any two numbers obeying
\begin{equation}\label{eRhos}
	0 < \rho < \rho' < \frac{1}{C}.
\end{equation}
Under these assumptions, there exists $h_0 > 0$ sufficiently small and $C_0 > 0$ sufficiently large for which, when $0 < h \leq h_0$ and $|z| \leq \rho,$
we have the following bound on the $L^2(\Bbb{R}^d)$ operator norm of the resolvent of $p^w(x,hD_x)$ at $z$:
\begin{equation}\label{eProp21Conclusion}
	||(p^w(x,hD_x) - z)^{-1}|| \leq e^{C_0h^{-d}} \mathop{\prod_{z_j\in \opnm{Spec}(p^w(x,hD_x))}}_{|z_j|\leq \rho'} |z-z_j|^{-1}.
\end{equation}
Furthermore, $\opnm{Spec}(p^w(x,hD_x)) \cap \{|z| \leq \rho'\}$ is a discrete set, and counting with multiplicity,
$$
	\#\left[(\opnm{Spec} p^w(x,hD_x))\cap\{|z|\leq \rho'\}\right] = \BigO(h^{-d}).
$$
\end{proposition}

\begin{proof}
Let $\chi \in C_0^\infty(\Bbb{R}^{2d})$ be a standard cutoff function adapted to $\{|(x,\xi)|\leq \tilde{C}\}$ for $\tilde{C}$ from (\ref{eEllinfty}), meaning that $\chi$ takes values in $[0,1]$ and
$\chi|_{\{|(x,\xi)|\leq \tilde{C}\}} \equiv 1.$ Let us define
$$
	\tilde{p}(x,\xi) = p(x,\xi) + C_1 \chi(x,\xi),
$$
where $C_1>0$ is taken so large that, when $|z|\leq \rho'$, we have
\begin{equation}
\label{ell_global}
\jvRe\left(\tilde{p}(x,\xi) - z\right) \geq \frac{m(x,\xi)}{{\cal O}(1)} \quad\quad \forall (x,\xi) \in \Bbb{R}^{2d}.
\end{equation}
Here we have used (\ref{eEllinfty}) and the fact that $m\geq 1$.

Henceforth we shall use capital letters to denote the $h$-dependent Weyl quantizations of $p$ and $\tilde{p},$ writing
$$
	P = p^w(x, hD_x), \quad \tilde{P} = \tilde{p}^w(x,hD_x).
$$
It follows from (\ref{eEllinfty}) that, when equipped with the domain $H(m)$, the operators $P$ and $\tilde{P}$ become closed and densely defined
on $L^2(\Bbb{R}^d)$, and our operator norms will refer to this $L^2(\Bbb{R}^d)$ as well. Now an application of (\ref{ell_global}) shows that, for $h$ sufficiently small, the family $\tilde{P}-z$ is invertible and
\begin{equation}\label{eptbdd}
	||(\tilde{P} - z)^{-1}|| = \BigO(1)
\end{equation}
uniformly for $|z|\leq \rho'$.  This fact may be seen by noting that
$$
	q(x,\xi;z) := (\tilde{p}(x,\xi)-z)^{-1}
$$
varies in a bounded set in $S(m^{-1})$ when $|z| \leq \rho'$.  The pseudodifferential symbol calculus, described in Proposition 7.7,
Theorem 7.9, and Theorem 7.11 in \cite{DiSj}, gives that
$$
	q^w(x,hD_x;z)(\tilde{p}^w(x,hD_x) - z) = 1 + \BigO(h),
$$
where the remainder estimate in the right hand side refers to the operator norm on $L^2(\Bbb{R}^d).$  The right-hand side may therefore
be inverted by the Neumann series, and it is also true that $q^w(x,hD_x; z)$ is bounded as a map from $L^2(\Bbb{R}^d)$ to $H(m) \subseteq L^2(\Bbb{R}^d)$. This demonstrates the existence of a
uniformly bounded left inverse for $\tilde{P}-z$ when $|z|\leq \rho'.$  The same procedure, multiplying with $q^w(x,hD_x;z)$ on the right instead,
completes the proof of (\ref{eptbdd}).

We can apply similar reasoning to the operator $P-z$ for any $z$ such that $p(x,\xi) - z$ is elliptic on all of $\Bbb{R}^{2d}$.  To fix matters,
let us choose $z = -\rho/2.$  Since $\jvRe p \geq 0,$ we have uniform ellipticity for $p(x,\xi) + \rho/2$, and therefore we have
\begin{equation}\label{epbdd}
	||(P+\rho/2)^{-1}|| = \BigO(1),
\end{equation}
for $h$ sufficiently small.

The resolvents of $P$ and $\tilde{P}$ are closely related, where both exist. We introduce the operators
\begin{equation}\label{eK1}
	K_1(z) = (\tilde{P}-z)^{-1}(P-\tilde{P}),
\end{equation}
\begin{equation}\label{eK2}
	K_2(z) = (P-z)^{-1}(\tilde{P}-P).
\end{equation}
It is straightforward to check that
\begin{equation}\label{eResolvents1}
	(P-z)^{-1} = (1+K_1(z))^{-1}(\tilde{P}-z)^{-1},
\end{equation}
\begin{equation}\label{eResolvents2}
	(\tilde{P}-z)^{-1} = (1+K_2(z))^{-1}(P-z)^{-1}.
\end{equation}
The operators $K_1(z)$ and $K_2(z)$ are both trace class operators, as will be shown below.

We shall now recall some general operator-theoretic facts regarding trace class operators.  We shall refer to \cite{GoGoKr90},
but an interested reader may also refer to \cite{Si05} or \cite{GoKr}.  A compact operator $A$ acting on a Hilbert space is of trace class if
its singular values, the non-vanishing eigenvalues of $(A^*A)^{1/2}$, are summable.  Henceforth assume that $A$ is of trace class. Write the
non-vanishing eigenvalues of $A$ as a sequence $\{\lambda_j\}_{j=1}^\infty,$ repeated for algebraic multiplicity, arranged so that
$|\lambda_1|\geq |\lambda_2|\geq \ldots\,$, and write the singular values as $\{s_j\}_{j=1}^\infty,$ also repeated for multiplicity and arranged in
decreasing order. Then
$$
	||A||_{tr} = \sum_{j=1}^\infty s_j < \infty.
$$
Note that pre- or post-multiplying by bounded operators has a limited effect on the trace norm of $A$: if $B, C$ are bounded operators on the same Hilbert space, then
\begin{equation}\label{eProdTr}
	||BAC||_{tr} \leq ||B||\,||A||_{tr}\,||C||,
\end{equation}
(\cite{GoGoKr90}, Proposition VI.4.2).
The summability of $\{s_j\}_{j=1}^\infty,$ which implies summability of $\{|\lambda_j|\}_{j=1}^\infty$ by a classical theorem of H. Weyl \cite{Weyl49}, allows us to extend the definition of matrix determinant to $1+A$, writing
\begin{equation}\label{eFredDetDef}
	\det(1+A) = \prod_{j=1}^\infty (1+\lambda_j).
\end{equation}
Here $1$ denotes the identity operator; see \cite{GoGoKr90}, Theorem VII.6.1.

We have the inequalities (\cite{GoGoKr90}, Theorem VII.3.3),
\begin{equation}\label{eDetBddByTr}
	|\det(1+A)| \leq \prod_{j=1}^\infty (1+s_j) \leq \exp\left(||A||_{tr}\right).
\end{equation}
The product formula for determinants, from the same theorem and remarks immediately following, gives us the usual formula
\begin{equation}\label{eDetInv}
	\frac{1}{\det(1+A)} = \det\left[(1+A)^{-1}\right],
\end{equation}
where the inverse exists, which coincides with when the determinant is nonzero. It follows therefore from (\ref{eResolvents1}) that
$z_0$ with $|z_0|\leq \rho'$ is in the spectrum of $P$ precisely when the determinant of $1+K_1(z)$ vanishes at $z=z_0$. Furthermore, the algebraic
multiplicity of $z_0$ is equal to the order of vanishing of the determinant --- see, e.g., \cite{MarkusPencils}, Lemma 7.1, for this
essentially well known consequence of the argument principle.

Finally, the determinant allows us to control the operator norm of the inverse of $1+A,$ using \cite{GoGoKr90}, Theorem X.1.1:
\begin{equation}\label{eDetBddsResolvent}
	||(1+A)^{-1}|| \leq \frac{1}{|\det(1+A)|}\prod_{j=1}^\infty (1+s_j) \leq \frac{1}{|\det(1+A)|} \exp\left(||A||_{tr}\right).
\end{equation}
This is the primary tool we will use to bound the resolvent $(P-z)^{-1},$ using (\ref{eResolvents1}).

To bound $||P-\tilde{P}||_{tr},$ we refer to \cite{DiSj}, Theorem 9.4, which states that, for a symbol $a \in \mathcal{S}'(\Bbb{R}^{2d}),$ we
have the non-semiclassical ($h=1$) bound
\begin{equation}\label{eCutoffTr}
	||a^w(x,D_x)||_{tr} \leq C_d \sum_{(j,k) \in \Bbb{Z}^{2d}} ||\widehat{\psi_{j,k} a}||_{L^1(\Bbb{R}^{2d})}.
\end{equation}
The Hilbert space acted on by $a^w(x,D_x)$ is $L^2(\Bbb{R}^d),$ as usual.  The nonnegative $C_0^\infty(\Bbb{R}^{2d})$ functions $\psi_{j,k}(x,\xi) = \psi(x-j, \xi-k)$ are translates of one another chosen to form a partition of unity:
$$
	\sum_{(j,k)\in \Bbb{Z}^{2d}}\psi_{j,k}(x,\xi) = 1, \quad \forall(x,\xi)\in \Bbb{R}^{2d}.
$$

As alluded to on page 113 of \cite{DiSj}, when $a$ is sufficiently differentiable we have
\begin{equation}\label{ePartialsTr}
	\sum_{(j,k)\in \Bbb{Z}^{2d}}||\widehat{\psi_{j,k} a}||_{L^1(\Bbb{R}^{2d})} \leq C_\psi \sum_{|\alpha| \leq 2d+1}
||\partial_{x,\xi}^\alpha a||_{L^1(\Bbb{R}^{2d})}.
\end{equation}
This may be verified by standard integration by parts, combined with the fact that
$$
	\sum_{(j,k)\in \Bbb{Z}^{2d}} |\partial^\alpha_{x,\xi}\psi_{j,k}(x,\xi)| = \BigO_\alpha(1).
$$
The trace norm of a non-semiclassical pseudodifferential operator $a^w(x, D_x)$ may therefore be bounded by the right-hand side of (\ref{ePartialsTr}).

To apply this result to a semiclassical quantization, write $a_h(x,\xi) = a(x,h\xi)$ and note that
$$
	a^w(x,hD_x) = (a_h)^w(x,D_x).
$$
Changing variables then gives that
$$
	\sum_{|\alpha| + |\beta| \leq 2d+1} ||\partial_x^\alpha \partial_\xi^\beta a_h(x,\xi)||_{L^1} = \sum_{|\alpha| + |\beta| \leq 2d+1} h^{-d+|\beta|}||\partial_x^\alpha \partial_\xi^\beta a(x,\xi)||_{L^1}.
$$
In the case where the derivatives of $a$ up to order $2d+1$ are in $L^1(\Bbb{R}^{2d})$, which is certainly the case if $a \in C_0^\infty(\Bbb{R}^{2d}),$ we use (\ref{eCutoffTr}), (\ref{ePartialsTr}) to obtain
$$
	||a^w(x,hD_x)||_{tr} = \BigO_a(h^{-d}), \quad h\rightarrow 0^+.
$$

Recall that $p(x,\xi) - \tilde{p}(x,\xi) \in C_0^\infty(\Bbb{R}^{2d})$ and that $||(\tilde{P}-z)^{-1}|| = \BigO(1)$ for $h$ sufficiently small and $|z| \leq \rho'$. Therefore, as a consequence of (\ref{eProdTr}), (\ref{eptbdd}), and the definition (\ref{eK1}) of $K_1$,
\begin{equation}\label{eK1trbdd}
	||K_1(z)||_{tr} = \BigO(h^{-d}),
\end{equation}
uniformly in $h$ when $h$ is sufficiently small and in $z$ when $|z| \leq \rho'.$  Similarly, now using (\ref{epbdd}) and (\ref{eK2}),
$$
	||K_2(-\rho/2)||_{tr} = \BigO(h^{-d}).
$$

Let us introduce the holomorphic function,
$$
	D(z) = \det(1+K_1(z)),
$$
defined on $\{|z| \leq \rho\}$ for $h > 0$ sufficiently small.  We wish to bound $|D(z)|$ from below on $\{|z| \leq \rho\}$, at which point we will apply (\ref{eDetBddsResolvent}) to complete the proof.
Using (\ref{eDetBddByTr}) and (\ref{eK1trbdd}), we see that
\begin{equation}\label{eDUppBd}
	|D(z)| \leq e^{C h^{-d}},
\end{equation}
when $h$ is sufficiently small, uniformly in $z$ when $|z| \leq \rho'.$  Similarly,
$$
	|\det(1+K_2(-\rho/2))| \leq e^{C h^{-d}}.
$$
It is straightforward to check, using (\ref{eResolvents1}) and (\ref{eResolvents2}), that when both $K_1$ and $K_2$ are defined,
$$
	(1+K_2(z))^{-1} = 1+K_1(z).
$$
Using the standard determinant formula for inverses (\ref{eDetInv}), we therefore have
\begin{equation}\label{eDLwrBd}
	|D(-\rho/2)| \geq e^{-Ch^{-d}}.
\end{equation}

We extend the lower bounds for $|D(z)|$ at $z = -\rho/2$ to the rest of $\{|z|\leq \rho\}$ using the standard complex analytic
Lemma \ref{BoundsOnHolomorphicFunctionsOnDisc} below.  To biholomorphically map $\{|z| \leq 1\}$ onto $\{|z|\leq \rho'\},$ we use
a M\"{o}bius transformation
$$
	\varphi(z) = \rho' \cdot \frac{2\rho' z - \rho}{2 \rho' - \rho z},
$$
chosen so that $\varphi(0) = -\rho/2.$  In this way, the function
$$
	u(z) = \frac{1}{D(-\rho/2)} D(\varphi(z))
$$
is holomorphic on $\{|z| < 1\}$ and continuous on $\{|z| \leq 1\}.$  Furthermore, clearly $u(0) = 1,$ and from (\ref{eDUppBd}) and (\ref{eDLwrBd})
we obtain
\begin{equation}\label{euUppBd}
	|u(z)| \leq e^{Ch^{-d}}, \quad \forall |z| \leq 1
\end{equation}
for $h$ sufficiently small.  Define
$$
	r = \max_{|z| = \rho} |\varphi^{-1}(z)|,
$$
noting that $r < 1$ as $\varphi^{-1}$ is continuous and $\rho < \rho'.$  Furthermore, since $\varphi$ depends only on $\rho$ and $\rho',$ the same holds true for $r.$  Also define
$$
	R = \frac{r+1}{2}.
$$
These choices are made to ensure that
\begin{equation}\label{evarphiInclusions}
	\{|\varphi(z)| \leq \rho\} \subseteq \{|z| \leq r\} \subseteq \{|z| \leq R\} \subseteq \{|\varphi(z)|\leq \rho'\}.
\end{equation}

We are in a position to apply Lemma \ref{BoundsOnHolomorphicFunctionsOnDisc} to $u(z)$ with $0 < r < R < 1$ and $\lambda = Ch^{-d}.$  Recall that the constants in (\ref{eDUppBd}) and (\ref{eDLwrBd}) depend only on ellipticity of $\tilde{p}-z$ taken over $\{|z| \leq \rho'\}$, ellipticity of $p+\rho/2,$ and bounds for $||P-\tilde{P}||_{tr}.$  Therefore the constant in $\lambda = Ch^{-d}$ in (\ref{euUppBd}) may be taken to depend only on $p$, $\rho$, and $\rho'.$  From this, we conclude that there exist constants depending only on $p, \rho,$ and $\rho'$ where, for $h$ sufficiently small,
\begin{equation}\label{eZerosOfDet}
	\#\{z \::\: D(\varphi(z)) = 0, |\varphi(z)| \leq \rho \} \leq Ch^{-d}
\end{equation}
(counting for multiplicity) and, for all $|z| \leq r$,
$$
	\frac{1}{|D(-\rho/2)|}|D(\varphi(z))| \geq e^{-Ch^{-d}} \mathop{\prod_{D(\varphi(z_j)) = 0}}_{|z_j|\leq R} |z - z_j|.
$$
Multiplying both sides by $|D(-\rho/2)|$ and recalling (\ref{eDLwrBd}) gives the slightly simpler formula
\begin{equation}\label{eLwrBdWPhi}
	|D(\varphi(z))| \geq e^{-Ch^{-d}} \mathop{\prod_{D(\varphi(z_j)) = 0}}_{|z_j|\leq R} |z - z_j|, \quad \forall |z|\leq r.
\end{equation}

As already observed, counting for multiplicity, we have
\begin{equation}\label{ePolesOfDAreSpec}
	\{\zeta \::\: D(\zeta) = 0, |\zeta| \leq \rho' \} = \opnm{Spec}(P)\cap \{|\zeta|\leq \rho'\}.
\end{equation}
Note that the reasoning leading up to (\ref{eZerosOfDet}) implies that
\begin{equation}\label{eZerosOfDetGeneral}
	\#\{\zeta \::\: D(\zeta) = 0, |\zeta| \leq \rho'\} = \BigO_{p, \rho'}(h^{-d}),
\end{equation}
simply by repeating the same argument with $0 < \rho < \rho' < 1$ replaced by $0 < \rho' < (\rho'+1)/2 < 1.$

Using $\varphi$ as a change of variables and recalling (\ref{evarphiInclusions}), we obtain from (\ref{eLwrBdWPhi}) the inequality
\begin{equation}\label{eLwrBdWOPhi}
	|D(z)| \geq e^{-Ch^{-d}} \mathop{\prod_{D(z_j) = 0}}_{|\varphi^{-1}(z_j)|\leq R} |\varphi^{-1}(z) - \varphi^{-1}(z_j)|, \quad \forall |z| \leq \rho.
\end{equation}
Also from (\ref{evarphiInclusions}), we have
$$
	\{z \::\: |z| \leq R\} \subseteq \{\varphi^{-1}(z) \::\: |z| \leq \rho'\}.
$$
Focusing on the product in (\ref{eLwrBdWOPhi}), we wish to include all zeros up to modulus $\rho'$ and eliminate $\varphi^{-1}.$  We are faced with the task of bounding from below
\begin{equation}\label{eExtendingProduct}
	\mathop{\prod_{D(z_j) = 0}}_{|\varphi^{-1}(z_j)|\leq R} |\varphi^{-1}(z) - \varphi^{-1}(z_j)|
	= \frac{\prod_{D(z_j) = 0,|z_j| \leq \rho'} |\varphi^{-1}(z) - \varphi^{-1}(z_j)|}
	{\prod_{D(z_j) = 0, |z_j| \leq \rho', |\varphi^{-1}(z_j)| > R} |\varphi^{-1}(z) - \varphi^{-1}(z_j)|}
\end{equation}
where $|z| \leq \rho.$  Since $\varphi^{-1}$ maps to the unit disc, each term in the denominator cannot exceed 2.  Furthermore, by (\ref{eZerosOfDetGeneral}), there are at most $\BigO(h^{-d})$ terms and therefore the denominator is bounded above:
\begin{equation}\label{eTwoBound}
	\mathop{\prod_{D(z_j) = 0}}_{|z_j|\leq \rho', |\varphi^{-1}(z_j)| > R}|\varphi^{-1}(z) - \varphi^{-1}(z_j)| \leq 2^{Ch^{-d}}.
\end{equation}
The mean value theorem gives that
$$
	|\varphi(\zeta) - \varphi(\zeta_j)| \leq ||\varphi'||_{L^\infty(\{|z|\leq 1\})}|\zeta - \zeta_j|,
$$
so
$$
	|z - z_j| \leq C_{\varphi}|\varphi^{-1}(z) - \varphi^{-1}(z_j)|,\quad C_{\varphi}>0.
$$
This may be applied to the numerator of (\ref{eExtendingProduct}), which also has at most $\BigO(h^{-d})$ terms.  Combined with (\ref{eTwoBound}), this gives
$$
	\mathop{\prod_{D(z_j) = 0}}_{|\varphi^{-1}(z_j)|\leq R} |\varphi^{-1}(z) - \varphi^{-1}(z_j)| \geq (2C_{\varphi})^{-Ch^{-d}} \mathop{\prod_{D(z_j) = 0}}_{|z_j| \leq \rho'} |z - z_j|
$$
Inserting this into (\ref{eLwrBdWOPhi}) and using the resolvent bound (\ref{eDetBddsResolvent}) gives
$$
	||(1+K_1(z))^{-1}|| \leq e^{Ch^{-d}}e^{||K_1(z)||_{tr}}\mathop{\prod_{D(z_j) = 0}}_{|z_j| \leq \rho'} |z - z_j|^{-1}.
$$
The conclusion (\ref{eProp21Conclusion}) of the proposition follows immediately from (\ref{eResolvents1}), (\ref{eptbdd}), (\ref{eK1trbdd}), and (\ref{ePolesOfDAreSpec}).

\end{proof}

We complete the proof by supplying the missing lemma referenced in the proof of Proposition \ref{ResolventUpperBounds}. See
\cite{MarkusPencils}, Chapter 1, and \cite{Sj01} for closely related arguments.

\begin{lemma}\label{BoundsOnHolomorphicFunctionsOnDisc}
Let $u(z):\{|z| \leq 1\}\rightarrow \Bbb{C}$ be continuous and holomorphic on $\{|z| < 1\}.$  Furthermore assume that $u(0) = 1$ and that $|u(z)| \leq e^{\lambda}$ when $|z| \leq 1,$ where $\lambda \geq 0.$  Let $Z(\rho)$ denote the set of zeros of $u$ for which $|z| < \rho,$ repeated according to multiplicity.  Then, for any $0 < r < R < 1,$ there exist constants depending only on $r$ and $R$ for which
\begin{equation}\label{eLem2.2p1}
    \#(Z(r)) \leq C\lambda
\end{equation}
and, for all $|z| < r,$
\begin{equation}\label{eLem2.2p2}
    |u(z)| \geq e^{-C\lambda}\prod_{z_j \in Z(R)}|z-z_j|.
\end{equation}
\end{lemma}

Note that the fact that $u$ is holomorphic and not identically zero immediately implies that $Z(R')$ is finite for any $0 < R' < 1.$

\begin{proof}
From Jensen's formula (see \cite{Ahlfors}, p. 207-208),
$$
    0 = \log|u(0)| = -\sum_{z_j \in Z(R)} \log\left(\frac{R}{|z_j|}\right) + \frac{1}{2\pi}\int_0^{2\pi} \log|u(Re^{i\theta})|\,d\theta.
$$
Since $\log(R/|z_j|) \geq 0$ for all $z_j \in Z(R)$ and since $Z(r) \subseteq Z(R),$
\begin{eqnarray*}
    \frac{1}{2\pi}\int_0^{2\pi} \log|u(Re^{i\theta})|\,d\theta
    &=& \sum_{z_j \in Z(R)}\log\left(\frac{R}{|z_j|}\right)
    \\ & \geq & \sum_{z_j\in Z(r)}\log \left(\frac{R}{|z_j|}\right)
    \\ &\geq & \log\left(\frac{R}{r}\right)(\#Z(r)).
\end{eqnarray*}
By our hypotheses on $u,$ we know that $\log|u(Re^{i\theta})| \leq \lambda,$ and so we obtain (\ref{eLem2.2p1}) with constant $C = (\log(R/r))^{-1}.$

Note that the constant in (\ref{eLem2.2p1}) may be taken to depend on $r$ alone, upon substituting $R = (1+r)/2.$

If $r', x \in [0, R],$ a change of variables shows that
$$
	-\int_{r'}^R \log|\rho-x|\,d\rho \leq -\int_{-1}^1 \log|\rho|\,d\rho = 2.
$$
From this, choosing $r' = \frac{1}{2}(R+r),$ we see that
$$
	-\int_{r'}^R \sum_{z_j \in Z(R)}\log\left|\rho - |z_j|\right|\,d\rho \leq 2(\#Z(R)).
$$
Because the integrand cannot be everywhere greater than the average value, there exists some $R' \in (r', R)$ where
$$
	-\sum_{z_j \in Z(R)} \log\left| R' - |z_j|\right| \leq \frac{2}{R-r'} (\#Z(R)) = \frac{4}{R-r}(\#Z(R)).
$$
Using (\ref{eLem2.2p1}) to bound $\#Z(R),$ we conclude that, for any $z$ with $|z| = R',$ we have
\begin{equation}\label{eLowerProdBd}
	\prod_{z_j \in Z(R)} |z-z_j| \geq \prod_{z_j \in Z(R)}\left| |z| - |z_j| \right| \geq e^{-\frac{4}{R-r}(\#Z(R))} \geq e^{-C\lambda}.
\end{equation}

Write $f(z)$ for the part of $u(z)$ with the zeros up to modulus $R$ removed,
\begin{equation}\label{eLem2.2Deff}
	f(z) = \frac{u(z)}{\prod_{z_j \in Z(R)} (z-z_j)}.
\end{equation}
Combining (\ref{eLowerProdBd}) with $|u(z)| \leq e^{\lambda}$ gives
$$
	|f(z)| \leq e^{C_0 \lambda}
$$
for some $C_0$ depending only on $r, R.$  This inequality holds whenever $|z| \leq R'$ by the maximum principle.  Therefore the harmonic function
\begin{equation}\label{eLem2.2DefG}
	G(z) = C_0\lambda - \log |f(z)|
\end{equation}
is nonnegative on $\{|z| \leq R' \}.$  Furthermore, since $u(0) = 1$ and $R < 1,$
$$
	G(0) = C_0\lambda + \sum_{z_j \in Z(R)} \log |z_j| \leq C_0\lambda.
$$

Applying Harnack's inequality (see \cite{Ahlfors}, p. 243) gives, for $|z| \leq r,$
$$
	G(z) \leq \frac{R'+r}{R'-r} G(0) \leq 2\frac{R+r}{R-r}C_0 \lambda.
$$
From this and the definition (\ref{eLem2.2DefG}) of $G$, we obtain that
$$
	\log |f(z)| \geq C_0(1-2\frac{R+r}{R-r})\lambda.
$$
Taking the exponential function of both sides and using the definition (\ref{eLem2.2Deff}) of $f$ proves the final statement in the lemma, (\ref{eLem2.2p2}).

\end{proof}

\section{Resolvent Estimates in the Quadratic Case}\label{SectionQuadraticCase}

The largest term contributing to the upper bound in Theorem \ref{MainTheorem} is the term coming from the quadratic part $q$ of the symbol $p.$
Therefore, to begin, we use a simple rescaling argument to find out for which $f(h)$ the resolvent estimate in Proposition
\ref{ResolventUpperBounds} gives a resolvent upper bound of $\BigO(h^{-1-\gamma})$ for $|z| \leq hf(h)$ when $z$ is sufficiently far from the spectrum of $q^w(x,hD_x).$

Conditions on distance from the spectrum are also delicate, as if no $|z| \leq hf(h)$ is sufficiently far from the spectrum, Theorem \ref{MainTheorem} is vacuous.  With this in mind, we find a satisfactory minimum value of $\operatorname{dist}(z,\operatorname{Spec} (q^w(x,hD_x))$ for which the
subset of $\{|z| \leq hf(h)\}$ sufficiently far from $\operatorname{Spec} (q^w(x,hD_x))$ has measure approaching the full measure of the set
$\{|z| \leq hf(h)\}$ in the limit $h \rightarrow 0^+.$

\begin{proposition}\label{UnconjugatedQuadratic}

Let $q(X):\Bbb{R}^{2d}\rightarrow \Bbb{C}$ be a quadratic form
obeying $\jvRe q(X) \geq |X|^2/M$ for some $M > 0$.  For any $\gamma
> 0$ there exists $C = C_\gamma>0$ such that, defining the slowly growing parameter $f(h)$ via $$f(h) =
\frac{1}{C}\left(\frac{\log\frac{1}{h}}{\log\log\frac{1}{h}}\right)^{1/d},$$
there exists $h_0 > 0$ where, for all $h \in (0,h_0]$ and $|z|
\leq hf(h)$ with $$\operatorname{dist}(z,
\operatorname{Spec}(q^w(x,hD_x)))\geq hf(h)^{(1-d)/2},$$ we have
the resolvent estimate
$$
	||u||\leq \BigO(h^{-1-\gamma}) ||(q^w(x,hD_x)-z)u||
$$ for any $u$ in the domain of $q^w(x,hD_x).$  Here, norms are taken in $L^2(\Bbb{R}^d).$

\end{proposition}

\begin{remark}\label{remGammaDependence}

The constant $C = C_\gamma$ in the definition of $f(h)$ has a simple dependence on $\gamma,$ and may be chosen $$C_\gamma =
\frac{C'}{\gamma^{1/d}},$$ with $C'$ coming from the constants in Proposition \ref{ResolventUpperBounds}.  Restrictions on values of $\gamma,$
like those in the main theorem, thus prevent us from expanding our region $\{|z| \leq hf(h)\}$ insofar as we cannot multiply by arbitrary
constants.

The domain of $q^w$ is $$\operatorname{Dom}(q^w(x,hD_x)) = \{u \in L^2(\Bbb{R}^d)\::\: (x^2 + (hD_x)^2)u \in L^2(\Bbb{R}^d)\},$$ and, as is shown
in Appendix \ref{WeightsAppendix}, for $u \in \operatorname{Dom}(q^w)$ the harmonic oscillator $x^2 + (hD_x)^2$ gives the same graph norm as
$q^w(x,hD_x)$ which may be any elliptic quadratic form: $$||u||_{L^2} + \frac{1}{h}||(x^2+(hD_x)^2)u||_{L^2} \sim ||u||_{L^2} + \frac{1}{h}||q^w(x,hD_x)u||_{L^2},$$
uniformly in $h$.

\end{remark}

\begin{proof}By rescaling, we seek $f(h)$ and a condition on distance from
the spectrum such that $$||(q^w(x,hD_x) - z)^{-1}||_{L^2(\Bbb{R}^d) \rightarrow L^2(\Bbb{R}^d)} \leq \BigO(h^{-1-\gamma}),$$ for $\gamma
> 0$ fixed, $|z| \leq hf(h)$, and $h$ sufficiently small.  Our primary tool will be
the upper bound in Proposition \ref{ResolventUpperBounds}.

In order to scale $z$ into a $h$-independent set, factor out
$hf(h),$ giving $$q^w(x,hD_x) - z =
hf(h)\left(\frac{1}{hf(h)}q^w(x,hD_x) - \frac{z}{hf(h)}\right).$$

Next, change variables as in (\ref{QuadraticChangeVariables}) with $\alpha = hf(h)$ and $\tilde{h} = h/\alpha = 1/f(h).$  Writing $\tilde{z} =
z/hf(h),$ the unitary equivalence in (\ref{QuadraticChangeVariables}) gives, with $L^2(\Bbb{R}^d)$ operator norms,
$$
	||(q^w(x,hD_x) - z)^{-1}|| = \frac{1}{hf(h)}||(q^w(x,\tilde{h}D_x) - \tilde{z})^{-1}||.
$$

Since we will have $f(h) \gg 1$ as $h \rightarrow 0,$ it is sufficient to choose $f(h)$ to establish
$$
	||(q^w(x,\tilde{h}D_x) - \tilde{z})^{-1}|| \leq \BigO(h^{-\gamma}).
$$

The dominant error in Proposition \ref{ResolventUpperBounds} will come from the product of distances from eigenvalues. For convenience, write
$$
	S_h(R) := \opnm{Spec}(q^w(x,hD_x)) \cap \{|z| \leq R\}.
$$
We will be interested in $S_{\tilde{h}}(1)$, and different regions may be related via the rule $S_h(R) = hS_1(h^{-1}R)$, which follows from (\ref{QuadraticChangeVariables}).  From Proposition \ref{ResolventUpperBounds}, we know that
$$
	\#(S_{\tilde{h}}(1)) = \BigO(\tilde{h}^{-d}), \quad\quad \tilde{h}\rightarrow 0^+,
$$
and this also may be verified from (\ref{eSjExplicitSpec}) as a consequence of the volume of a $d$-dimensional simplex.

If $\min_{\tilde{z}_j\in S_{\tilde{h}}(1)} |\tilde{z} - \tilde{z}_j| \geq D(h),$ we therefore have
$$
	\prod_{\tilde{z}_j \in S_{\tilde{h}}(1)} |\tilde{z} - \tilde{z}_j|^{-1} \leq \left(\frac{1}{D(h)}\right)^{C_2 \tilde{h}^{-d}}.
$$
On the other hand, the volume in $\Bbb{C}$ of the set
$$
	\left\{\tilde{z}\::\:\min_{\tilde{z_j}\in S_{\tilde{h}}(1)} |\tilde{z} - \tilde{z}_j| < D(h)\right\}
	= \bigcup_{\tilde{z}_j\in S_{\tilde{h}}(1)} B(\tilde{z}_j; D(h))
$$
is at most $\BigO(\tilde{h}^{-d})\pi D(h)^2.$  To have this volume tend to zero, any power larger than $d/2$ of $\tilde{h}$ will suffice for $D(h),$ but let us say that
$$
	D(h) = \tilde{h}^{(d+1)/2} = f(h)^{-(d+1)/2}.
$$

Having established a condition in terms of $\tilde{h}$ after rescaling, we undo the rescaling to obtain the corresponding set where $h$ is the semiclassical parameter in $S_h(hf(h))$:
\begin{multline}\label{eScalingExceptedSet}
	\left\{\tilde{z}\::\:|\tilde{z}|\leq 1,~\min_{\tilde{z_j}\in S_{\tilde{h}}(1)} |\tilde{z} - \tilde{z}_j| < f(h)^{-(d+1)/2}\right\} 
	\\ = \frac{1}{hf(h)}\left\{z\::\: |z|\leq hf(h),~ \min_{z_j\in S_h(hf(h))} |z-z_j| < hf(h)^{(1-d)/2}\right\}.
\end{multline}
Since the volume of the set on the left is $o(1)$ as $h\rightarrow 0^+,$ the volume on the set on the right is a vanishingly small fraction of the volume of $\{|z| \leq hf(h)\}.$  Consequently, the proportion of $z$ in $\{|z| \leq hf(h)\}$ which satisfy the hypotheses of Proposition \ref{UnconjugatedQuadratic}, and also of Theorem \ref{MainTheorem}, approaches $1$ as $h\rightarrow 0^+.$

Our choice of $D(h)$ gives the following upper bound for the product of distances from eigenvalues: when $|\tilde{z}| \leq 1,$
$$
	\prod_{\tilde{z}_j \in S_{\tilde{h}}(1)} |\tilde{z}-\tilde{z}_j|^{-1} \leq f(h)^{\frac{d+1}{2}\BigO(f(h)^d)} \leq e^{C_2 f(h)^d\log f(h)}.
$$
We combine this with Proposition \ref{ResolventUpperBounds}, where we are free to choose $\rho = 1$ and $\rho' = 2$.  This follows from the hypothesis
$$
	\jvRe q(x,\xi) \geq \frac{1}{M}|(x,\xi)|^2,
$$
from which we make the choice of order function
$$
	m(x,\xi) = 1 + \frac{1}{4M}|(x,\xi)|^2.
$$
This assures us that $\jvRe q(x,\xi) \geq 3m(x,\xi)$ for $|(x,\xi)|$ sufficiently large.  Putting this into Proposition \ref{ResolventUpperBounds} gives the upper bound
$$
	||(q^w(x,\tilde{h}D_x) - \tilde{z})^{-1}|| \leq \frac{1}{C_0}e^{C_1 f(h)^d + C_2 f(h)^d\log f(h)},
$$
valid for all $\tilde{z}$ where $|\tilde{z}| \leq 1$.  Since $f(h) \rightarrow \infty,$ the second term in the exponent dominates, and so it is sufficient to choose $f(h)$ where
$$
	e^{C_3 f(h)^d\log f(h)} \leq h^{-\gamma}
$$
for $h$ sufficiently small.

Taking logarithms gives
\begin{equation}\label{Findingf1}
	f(h)^d\log f(h) \leq \frac{\gamma}{C_3} \log\frac{1}{h}.
\end{equation}
The left-hand side is increasing in $f(h),$ and so there is a unique solution to where the inequality becomes an equality.  Write $F(h)$ for this solution, and note that $F(h)$ must increase as $h \rightarrow 0^+$ since the right-hand side of (\ref{Findingf1}) increases in the same limit.  We are able to determine the growth of $F(h)$ up to a constant factor, and for simpler formulas, we define $f(h)$ by taking a lower bound.

Since $F(h) \gg 1,$ we have $\log F(h) \in (1, F(h)^{\varepsilon})$ for any $\varepsilon > 0$ when $F(h)$ becomes sufficiently large depending on $\eps.$  Plugging in to (\ref{Findingf1}) gives
$$
	\left(\frac{\gamma}{C_3}\log\frac{1}{h}\right)^{1/(d+\varepsilon)} \leq F(h) \leq \left(\frac{\gamma}{C_3}\log\frac{1}{h}\right)^{1/d}.
$$
Taking logarithms of the above equation and noting that $\log\frac{\gamma}{C_3} \ll \log\log\frac{1}{h}$ gives that, for some $h_0 > 0,$ when $h \in (0, h_0]$ we have
$$
	\frac{1}{d+2\eps}\log\log \frac{1}{h} \leq \log F(h) \leq \frac{1}{d-\eps}\log\log\frac{1}{h}.
$$
Inserting into (\ref{Findingf1}) and solving gives
$$
	F(h) = \left(\frac{\gamma}{\tilde{C}(h)}\frac{\log\frac{1}{h}}{\log\log\frac{1}{h}}\right)^{1/d},
$$
for
$$
	\tilde{C}(h) \in \left(\frac{C_3}{d+2\varepsilon},\frac{C_3}{d-\varepsilon}\right)
$$
when $\varepsilon > 0$ and $h$ sufficiently small.

Therefore, up to a constant, our choice of $f(h)$ as in the proposition and in (\ref{Definingf}) is the only possible choice, given this
rescaling argument. To fix the idea, substitute $\varepsilon = 1/2$ and so define
$$
	C_\gamma = \left(\frac{C_3}{\gamma(d - 1/2)}\right)^{1/d}
$$
in the definition (\ref{Definingf}),
$$
	f(h) = \frac{1}{C_\gamma}\left(\frac{\log\frac{1}{h}}{\log\log\frac{1}{h}}\right)^{1/d}.
$$

\end{proof}

\subsection{The FBI-Bargmann Side}

It is convenient to restate the result in terms of the FBI-Bargmann transform, a unitary isomorphism
$$
	T = T_h:L^2(\Bbb{R}^d)\rightarrow H_{\Phi_0}(\Bbb{C}^d; h) = \operatorname{Hol}(\Bbb{C}^d) \cap L^2_{\Phi_0}(\Bbb{C}^d;h)
$$
where $$L^2_{\Phi_0}(\Bbb{C}^d;h) =
L^2(\Bbb{C}^d, e^{-2\Phi_0/h}dL(x)).$$  Recall that $dL(x)$ refers to Lebesgue measure on $\Bbb{C}^d.$  We refer the reader to
Appendix \ref{ChangeVariablesAppendix} for details.

It is sufficient for the current application to use the phase $\varphi_0(x,y) = \frac{i}{2}(x-y)^2,$ which gives the weight function $\Phi_0(x) = \frac{1}{2}(\jvIm x)^2.$  To
distinguish the symbols on the FBI-Bargmann side we will use fraktur letters, so for the full symbol we will write $$p^w(x,hD_x) = T^{-1}\frakp^w(x,hD_x)T$$ and for the quadratic part we will write
$$q^w(x,hD_x) = T^{-1}\frakq^w(x,hD_x)T.$$  The symbols on the FBI-Bargmann side are given by the exact Egorov relation 
$$
	\frakp(x,\xi) = p\circ\kappa_T^{-1}(x,\xi)
$$
for $\kappa_T(x,\xi) = (x-i\xi, \xi).$

Because $\kappa_T$ is linear, the quadratic part of our symbol on the FBI-Bargmann side, $\frakq(x,\xi),$ remains quadratic.  The contour for which the Weyl
quantization on the FBI-Bargmann side is defined may be shifted as in (\ref{ContourRealizationOfWeyl}), in which case we use almost analytic
extensions of $\frakp$ off $\Lambda_{\Phi_0} = \kappa_T(\Bbb{R}^{2d}),$ while $\frakq$ is already analytic.  Details regarding such extensions may be found, among other places, in Chapter 8 of M. Dimassi and
J. Sj\"{o}strand \cite{DiSj}.

Because $T$ is unitary, the FBI-Bargmann version of Proposition \ref{UnconjugatedQuadratic} is immediate.

\begin{proposition}\label{FBIQuadraticEstimate}

Let $T$ be the FBI-Bargmann transform defined for $\varphi_0$ in Appendix \ref{ChangeVariablesAppendix}, and let $\frakq(x,\xi) = q(\kappa_T^{-1}(x,\xi))$ be the quadratic form on $\Lambda_{\Phi_0} \subseteq \Bbb{C}^{2d}$ corresponding to $q(x,\xi)$.  For $|z| \leq hf(h)$ and $\operatorname{dist}(z,\operatorname{Spec}(\frakq^w))\geq hf(h)^{(1-d)/2},$ with $f(h)$ defined by (\ref{Definingf}), there
exists $h_0 > 0$ such that, when $h \in (0, h_0],$ $$||u|| \leq \BigO(h^{-1-\gamma})||(\frakq^w(x,hD_x) - z)u||,$$ for any $u\in
\operatorname{Dom}(\frakq^w)\subset H_{\Phi_0}(\Bbb{C}^d;h)$.  Here, norms are in $H_{\Phi_0}(\Bbb{C}^d;h)$.

\end{proposition}

\begin{remark}

As shown in Appendix \ref{WeightsAppendix}, we have the characterization $$\operatorname{Dom}(\frakq^w(x,hD_x)) = \{u \in H_{\Phi_0}(\Bbb{C}^d;h)\::\:
(1+|x|^2)u \in L^2_{\Phi_0}(\Bbb{C}^d;h)\}.$$

\end{remark}

\section{Localized Resolvent Estimates in the Quadratic Case}\label{SectionLocalizing}

The quadratic case only relates to the full symbol in a neighborhood of the origin, and so we introduce a cutoff function localizing near $0 \in \Bbb{C}^d$ and bound
the error.  We follow the proof of Proposition 5.2 in \cite{HeSjSt05}, and incorporate by reference two facts proven for elliptic quadratic differential
operators on the FBI-Bargmann side, using $\Pi$ for orthogonal projection from $L^2_{\Phi_0}(\Bbb{C}^d;h)$ onto $H_{\Phi_0}(\Bbb{C}^d;h).$  Equations (5.24) and (5.29) in
\cite{HeSjSt05} with $k = 0$ give us, with norms in $L^2_{\Phi_0}(\Bbb{C}^d;h)$,
\begin{equation}\label{CutoffAntiholomorphic}
	||(1-\Pi)(\chi_0 u)||\leq \BigO(h^{1/2})||u\bar{\partial}\chi_0||
\end{equation} and
\begin{equation}\label{CutoffCommutingWithPi}
	||[\frakq^w(x,hD_x),\Pi \chi_0]u|| \leq \BigO(h)||\un_K u||,
\end{equation} for any $u \in H_{\Phi_0}(\Bbb{C}^d;h)$, for any $\chi_0 \in C_0^\infty(\Bbb{C}^d)$, and for $K$ any fixed compact neighborhood of
$\operatorname{supp}\nabla \chi_0.$  (Constants depend on $\chi_0$, $K$, and $\frakq$.)  Here, $\un_K$ is the indicator function of the set $K$, and the operator $\Pi\chi_0$ in (\ref{CutoffCommutingWithPi}) refers to $u \mapsto \Pi(\chi_0 u)$.

\begin{proposition}\label{CutoffQuadraticEstimate}

Fix $\chi_0 \in C_0^\infty(\Bbb{C}^d)$ with $\chi_0 = 1$ near 0. As before, assume that $|z|\leq hf(h)$ and $$\operatorname{dist}(\operatorname{Spec}
(q^w(x,hD_x)), z) \geq hf(h)^{(1-d)/2},$$ for $f(h)$ defined in (\ref{Definingf}). Let $K$ be a compact neighborhood of $\operatorname{supp}\nabla\chi_0.$ Then, for
$\frakq^w = \frakq^w(x,hD_x)$, and for $h \in (0, h_0]$ where $h_0$ is taken sufficiently small,
$$
	||\chi_0 u || \leq \BigO(h^{-1-\gamma})||\chi_0(\frakq^w -z)u|| + \BigO(h^{-\gamma})||\un_K u||.
$$
Here $u \in H_{\Phi_0}(\Bbb{C}^d;h)$ and norms are taken in $L^2_{\Phi_0}(\Bbb{C}^d;h).$

\end{proposition}

\begin{remark}\label{remProjCutoffsDecay}

Note that, since Proposition \ref{FBIQuadraticEstimate} is applied to $\Pi(\chi_0 u)$ with $u \in H_{\Phi_0}(\Bbb{C}^d;h)$ and
$\chi_0 \in C^\infty_0(\Bbb{C}^d),$
it is necessary that $\Pi(\chi_0 u) \in \operatorname{Dom}(\frakq^w).$  We shall show a more general fact: writing
$$
	\Pi: L^2_\Phi(\Bbb{C}^d;h) \rightarrow H_\Phi(\Bbb{C}^d;h)
$$
for orthogonal projection, with a general $\Phi$ coming from a general quadratic phase $\varphi$ (see Appendix \ref{ChangeVariablesAppendix}), we
have that $\Pi$ maps any $\psi \in C_0^\infty(\Bbb{C}^d)$ to a holomorphic function with 
$$
	|x|^N\Pi\psi \in L^2_{\Phi}(\Bbb{C}^d;h), \quad \forall N \geq 0.
$$
Taking this with $\Phi = \Phi_0$ and the characterization of the domain of $\frakq(x,hD_x)$ given in
Appendix \ref{WeightsAppendix} shows that $\Pi(\chi_0 u) \in \opnm{Dom}(\frakq^w),$ which is used in the proof
of Proposition \ref{CutoffQuadraticEstimate} below.

For $\varphi$ a general quadratic phase and $\Phi$ the associated weight function we recall from Section 1 of \cite{Sj96} that
$$
	\Pi u(x) = Ch^{-d} \int_{\Bbb{C}^d} e^{\frac{2}{h}(\Psi(x,\bar{y})-\Phi(y))}u(y)\,dL(y),
$$
where $\Psi$ is a the
unique holomorphic quadratic form on $\Bbb{C}^{2d}$ fulfilling $$\Psi(x,\bar{x}) = \Phi(x)$$ for all $x \in \Bbb{C}^d.$  To bound an integral operator
with kernel $K(x,y)$ from $L^2_{\Phi}(\Bbb{C}^d;h)$ to itself it is sufficient to bound the integral operator with kernel
\begin{equation}\label{EqnConjugatedKernel}\tilde{K}(x,y) = e^{-\Phi(x)/h}K(x,y)e^{\Phi(y)/h}\end{equation} on $L^2(\Bbb{C}^d, dL),$ as in (\ref{eConjugatedKernel}) from Appendix \ref{SchurAppendix}.

The integral kernel for $|x|^N\Pi$ is thereore
$$
	\tilde{K}(x,y) = Ch^{-d}|x|^N \exp\left[\frac{1}{h}(-\Phi(x)+2\Psi(x,\bar{y}) - \Phi(y))\right].
$$
Since we are applying $|x|^N\Pi$ to $\psi,$ we can multiply $\tilde{K}(x,y)$ by $\un_{\opnm{supp}\psi}(y)$ without changing the result of $|x|^N\Pi\psi.$
Recalling also from \cite{Sj96}, Section 1, that
$$
	 - \Phi(x) + 2\jvRe \Psi(x,\bar{y}) - \Phi(y) \sim -|x-y|^2,
$$
we can apply Schur's test to bound an operator on $L^2(\Bbb{C}^d, dL)$ which has an integral kernel bounded by
$$
	|\tilde{K}(x,y)| \leq Ch^{-d}|x|^N e^{-\frac{1}{Ch}|x-y|^2}\un_{\operatorname{supp}\psi}(y).
$$

Since $|y| \leq \BigO(1)$ on $\operatorname{supp}\psi,$ in the integral $\langle x \rangle \sim \langle x-y\rangle,$ writing $\langle x \rangle =
(1+ |x|^2)^{1/2},$ and so changing variables gives
\begin{multline*}
	\int_{\Bbb{C}^d} |\tilde{K}(x,y)|\,(dL(x)\mathrm{~or~}dL(y))
	\\ \leq  \BigO(1)\int_{\Bbb{C}^d} h^{-d}\langle x-y\rangle^N e^{-\frac{1}{Ch}|x-y|^2}\,(dL(x)\mathrm{~or~}dL(y))
	\\ \leq \BigO(1)\int_{\Bbb{C}^d} h^{-d}\langle z \rangle^N e^{-\frac{1}{Ch}|z|^2}\,dL(z)  = \BigO(1).
\end{multline*}
We conclude that, in addition to $\Pi \psi \in H_{\Phi}(\Bbb{C}^d;h)$ by definition of $\Pi$, we furthermore have that
$$
	|x|^N\Pi\psi \in L^2_\Phi(\Bbb{C}^d;h), \quad \forall N \geq 0.
$$

As proven in Appendix \ref{WeightsAppendix}, $N=0$ and $N=2$ suffice to show $\Pi (\chi_0 u) \in \operatorname{Dom}(\frakq^w).$

\end{remark}

We now prove Proposition \ref{CutoffQuadraticEstimate}.

\begin{proof}

Using (\ref{CutoffAntiholomorphic}) and Proposition \ref{FBIQuadraticEstimate}, with norms in $L^2_{\Phi_0}(\Bbb{C}^d;h)$, \begin{eqnarray*}||\chi_0 u || &\leq& ||\Pi(\chi_0 u)|| +
||(1-\Pi)(\chi_0 u)|| \\ &\leq& \BigO(h^{-1-\gamma})||(\frakq^w -
z)\Pi\chi_0 u|| + \BigO(h^{1/2})||\un_K u|| \\ &\leq&
\BigO(h^{-1-\gamma})||\Pi\chi_0(\frakq^w - z)u|| +
\BigO(h^{-1-\gamma})||[\frakq^w, \Pi\chi_0]u|| + \BigO(h^{1/2})||\un_K
u||.\end{eqnarray*}

The proposition then immediately follows from (\ref{CutoffCommutingWithPi}) and the fact that $\Pi$ is an orthogonal projection.

\end{proof}

While the cutoff function in Proposition \ref{CutoffQuadraticEstimate} is fixed, the quadratic part of the symbol $\frakp$ is only a useful approximation
in a small region -- as we will see, a region of size $|x| \leq Ch^{\rho/2}$ for $\rho \in (2/3, 1).$  The change of variables framework
from Appendix \ref{ChangeVarsSubsection} will allow us to perform this rescaling in a systematic way.

We will use $\alpha = h^\rho$ and therefore $\tilde{h} = h/\alpha = h^{1-\rho}.$  To simplify notation, we denote our semiclassical dependence by writing $||\cdot||_h$ instead of $||\cdot||_{L^2_{\Phi_0}(\Bbb{C}^d;h)}.$  We apply Proposition \ref{CutoffQuadraticEstimate} to $\frakU_\alpha u = \alpha^{d/2}u(\alpha^{1/2}\cdot)$ in the
semiclassical regime with parameter $\tilde{h},$ and in doing so we make the assumption that $|\tilde{z}|\leq \tilde{h}f(\tilde{h}),$ where
$$
	f(\tilde{h}) = \frac{1}{C_\gamma}\left(\frac{\log\frac{1}{\tilde{h}}}{\log\log\frac{1}{\tilde{h}}}\right)^{1/d} = \frac{1}{C_\gamma}\left(\frac{(1-\rho)\log\frac{1}{h}}{\log(1-\rho) + \log\log\frac{1}{h}}\right)^{1/d}.
$$
Because $|\log(1-\rho)| \ll \log\log\frac{1}{h}$ for any fixed $\rho\in (0,1),$ we notice that
\begin{equation}\label{eftilde}
	f(\tilde{h}) = ((1-\rho)^{1/d}+o(1)) f(h)
\end{equation}
in the limit $h \rightarrow 0^+.$

Where $|\tilde{z}| \leq \tilde{h}f(\tilde{h})$ with $\operatorname{dist}(\operatorname{Spec}
(q^w(x,\tilde{h}D_x)), \tilde{z}) \geq \tilde{h}f(\tilde{h})^{(1-d)/2}$, the estimate in Proposition
\ref{CutoffQuadraticEstimate} then becomes
\begin{equation}\label{eRescaledCutoffIntermed}
||\chi_0\frakU_\alpha u||_{\tilde{h}} \leq
\BigO(\tilde{h}^{-1-\gamma})||\chi_0(\frakq^w(x,\tilde{h}D_x)-\tilde{z})\frakU_\alpha
u||_{\tilde{h}} + \BigO(\tilde{h}^{-\gamma})||\un_K\frakU_\alpha
u||_{\tilde{h}}.
\end{equation}
Using (\ref{QuadraticChangeVariables}) for quadratic differential operators, we obtain
\begin{equation}\label{eRescalingfrakq}
	\frakU_\alpha^{-1}(\frakq^w(x,\tilde{h}D_x)-\tilde{z})\frakU_\alpha u = \alpha^{-1}(\frakq^w(x, hD_x) - \alpha \tilde{z})u.
\end{equation}
Since $\frakU_\alpha^{-1}: H_\Phi(\Bbb{C}^d;\tilde{h}) \rightarrow H_\Phi(\Bbb{C}^d;h)$ is unitary, we may apply it to every term in (\ref{eRescaledCutoffIntermed}) and distribute over products using (\ref{ChangeVariablesDistributeProducts}).  Also using (\ref{eRescalingfrakq}), we arrive at
\begin{multline*}||\chi_0(\alpha^{-1/2}x)u||_h \leq \BigO(\alpha^{-1}\tilde{h}^{-1-\gamma})||\chi_0(\alpha^{-1/2}x)(\frakq^w(x,hD_x) - \alpha\tilde{z})u||_h \\ + \BigO(\tilde{h}^{-\gamma})||\un_K(\alpha^{-1/2}x)u||_h.
\end{multline*}
Recall $\tilde{h} = h/\alpha = h^{1-\rho},$
and so
$$
	\alpha^{-1}\tilde{h}^{-1-\gamma} = h^{-1-\gamma(1-\rho)}.
$$
Since we have assumed
$$
	|\tilde{z}| \leq \tilde{h}f(\tilde{h}),
$$
we see from (\ref{eftilde}) that, writing $z = \alpha \tilde{z}$,
$$
	|z| \leq hf(h)
$$
for $h$ sufficiently small and a change in the constant $C_\gamma$ in the definition (\ref{Definingf}) for $f$. Specifically, the new $C_\gamma$ may be defined by multiplying the old constant by $2(1-\rho)^{1/d}$.  The spectral avoidance condition $\operatorname{dist}(\tilde{z},\operatorname{Spec}(\frakq^w(x,\tilde{h}D_x))) \geq \tilde{h} f(\tilde{h})^{(1-d)/2},$
paired with the unitary equivalence between $\frakq^w(x,\tilde{h}D_x)$ and $\alpha^{-1}\frakq^w(x,hD_x),$ gives
$$
	\operatorname{dist}\left(z, \operatorname{Spec}(\frakq^w(x,hD_x))\right) \geq h f(h)^{(1-d)/2}
$$
provided a similar change of constant in the definition of $f$.  This completes the transition to the rescaled proposition which follows.

\begin{proposition}\label{RescaledCutoff}

Let $\chi_0 \in C_0^\infty(\Bbb{C}^d)$ be a cutoff function taking values in $[0,1]$ with $\chi(x) \equiv 1$ in a neighborhood of $0 \in \Bbb{C},$ and let $K$ be any compact neighborhood of $\opnm{supp}(\nabla \chi_0).$  Defining $f$ as in {\rm (\ref{Definingf})}, when $|z| \leq hf(h)$ with $\operatorname{dist}(z, \operatorname{Spec}(\frakq^w(x,hD_x)))\geq hf(h)^{(1-d)/2},$ we have the following rescaled resolvent estimate for $u \in H_{\Phi_0}(\Bbb{C}^d;h)$ and $h$ sufficiently small:
$$
	||\chi_0(h^{-\rho/2}x)u|| \leq \BigO(h^{-1-\gamma(1-\rho)})||\chi_0(h^{-\rho/2}x)(\frakq^w - z)u|| +
\BigO(h^{-\gamma(1-\rho)})||\un_K(h^{-\rho/2}x)u||.
$$
Here, norms are taken in $L^2_{\Phi_0}(\Bbb{C}^d;h).$

\end{proposition}

\begin{remark}\label{RemarkfDependenceGammaRho}

In view of (\ref{eftilde}), we see that the dependence of $f(h)$ on $\gamma(1-\rho)$ in Proposition \ref{RescaledCutoff} is the same as the dependence of $f(h)$ on
$\gamma$ in Proposition \ref{UnconjugatedQuadratic} (see Remark \ref{remGammaDependence}).  Specifically,
$$
	f(h) = \frac{(\gamma(1-\rho))^{1/d}}{C_0}\left(\frac{\log\frac{1}{h}}{\log\log\frac{1}{h}}\right)^{1/d},
$$
where $C_0$ does not depend on $\gamma$ or $\rho.$

\end{remark}

Finally, to apply the facts regarding the quadratic part of the symbol, $\frakq$, to the full symbol, $\frakp,$ we will need the following approximation lemma.

\begin{lemma}\label{LocalFullSymbol} Let $\frakp, \frakq$ be the FBI-side analogues, via the Egorov relation (\ref{eExactEgorov}), of symbols $p, q$ taken as in Theorem \ref{MainTheorem}. Because 
$$
p(x,\xi) - q(x, \xi) = \BigO(|(x,\xi)|^3),
$$
for $\varepsilon \geq h^{1/2}$ and $u \in H_{\Phi_0}(\Bbb{C}^d;h)$ we have
$$||\un_{\{|x| \leq \varepsilon\}}(\frakp^w - \frakq^w) u||_{L^2_{\Phi_0}(\Bbb{C}^d;h)} \leq \BigO(\varepsilon^3)||u||_{H_{\Phi_0}(\Bbb{C}^d;h)}.$$

\end{lemma}

\begin{proof}

Because $\kappa_T$ is a linear bijection from $\Bbb{R}^{2d}$ to $\Lambda_{\Phi_0}$, we certainly have that have that
\begin{equation}\label{ePQVanishOrder3}
	\mathfrak{p}(x,\xi) - \mathfrak{q}(x,\xi) = \BigO(|(x,\xi)|^3), \quad \forall (x,\xi) \in \Lambda_{\Phi_0}.
\end{equation}
Note that this estimate extends from $\Lambda_{\Phi_0}$ to all of $\Bbb{C}^{2d}$ when taking an almost analytic extension of $\frakp$, simply because an almost analytic extension of $\frakp-\frakq$ must have zero derivatives of first and second order at $(x,\xi) = (0,0).$
Furthermore,
\begin{equation}\label{eTriIneqOrder3}
	(|x+y|^2)^{3/2} \leq \BigO(1)(|x|^3 + |y|^3) \leq \BigO(1) (|x|^3 + |x-y|^3).
\end{equation}
Realize $\frakp^w - \frakq^w$ via the contour
$$
	(2\pi h)^{-d}\iint_{\Gamma_1(x)} e^{\frac{i}{h}(x-y)\cdot \theta} \psi_0(x-y)\left(\frakp(\frac{x+y}{2},\theta) - \frakq(\frac{x+y}{2}, \theta)\right)u(y)\,dy\,d\theta
$$
as in (\ref{ContourRealizationOfWeyl}).  Multiplying by $\un_{\{|x|\leq \eps\}},$ this can be extended to an integral kernel acting on $L^2(\Bbb{C}^d, dL)$ as in Appendix \ref{SchurAppendix}, giving a kernel of type (\ref{eSchurWeylKernel}) with
$$
	b(x,y,\theta) = \un_{\{|x|\leq \eps\}}(x)\psi_0(x-y)\left(\frakp(\frac{x+y}{2}, \theta) - \frakq(\frac{x+y}{2}, \theta)\right).
$$
We also note that, along $\Gamma_1$ defined in (\ref{eDefGammat}), linearity of $\partial_x \Phi_0$ gives
$$
	|(\frac{x+y}{2},\theta)| \leq \frac{|x+y|}{2} + ||\partial_x\Phi_0||\,|x+y| + |x-y| \leq \BigO(1)(|x+y| + |x-y|).
$$
Combining with (\ref{ePQVanishOrder3}) and (\ref{eTriIneqOrder3}), we see that
$$
	\left|\frakp(\frac{x+y}{2}, \theta) - \frakq(\frac{x+y}{2}, \theta)\right| \leq \BigO(1)(|x|^3 + |x-y|^3)
$$
along $\Gamma_1,$ and we therefore have an upper bound for our integral kernel $\tilde{K}(x,y)$ from (\ref{eSchurWeylFinalUpBd}),
$$
	|\tilde{K}(x,y)| \leq \BigO(1)h^{-d}\un_{\{|x|\leq \eps\}}(x)(|x|^3 + |x-y|^3)e^{-\frac{1}{h}|x-y|^2}.
$$

We then apply Schur's test to the two terms given by $|x|^3$ and $|x-y|^3.$  For $|x|^3,$ we use the bound $|x| \leq \eps,$ and by a change of variables
$$
	||\BigO(1)\eps^3 h^{-d}e^{-\frac{1}{h}|x-y|^2}||_{L^1(\Bbb{C}^d, dL)} = \BigO(\eps^3)
$$
regardless of whether we integrate in $x$ or in $y.$  For the term given by $|x-y|^3,$ we change variables to $z = (x-y)/h^{1/2}$ and obtain
$$
	||\BigO(1)|x-y|^3 h^{-d}e^{-\frac{1}{h}|x-y|^2}||_{L^1(\Bbb{C}^d, dL)} = \BigO(h^{3/2}),
$$
again regardless of which variable we integrate in.  By the hypothesis that $\eps \geq h^{1/2},$ the bound from Schur's test is $\BigO(\eps^3).$

\end{proof}

\section{Estimates in the Exterior Region}\label{SectionExterior}

In the region where $\chi_0(h^{-\rho/2}x) \neq 1,$ away from the origin, the symbol $\frakp(x,\xi)$ has real part bounded
 from below by $h^{\rho}/C.$  This allows us to use $\frakp^w(x,hD_x)$ to obtain upper bounds for portions of $u\in H_{\Phi_0}$ localized away from a neighborhood of the origin of size $h^{\rho/2}.$  The tool to obtain these bounds from lower bounds on $\frakp(x,\xi)$ is the quantization-multiplication Lemma \ref{QuantMult}.

\begin{proposition}\label{ExteriorEstimate} Let $\chi \in
C^\infty(\Bbb{C}^d)$ be fixed, taking on values in $[0,1],$ equal
to zero in a neighborhood of $0 \in \Bbb{C}^d,$ and equal to 1 for
$|x| > C$. Assume that $|z| \leq hf(h)$ for $f(h)$ defined in (\ref{Definingf}).  Then, for $\rho < 1$ and $u \in H_{\Phi_0}(\Bbb{C}^d;h),$ we have
\begin{equation}\label{ExteriorEstimateEqn}
	\int_{\Bbb{C}^d}\chi(h^{-\rho/2}x)|u|^2e^{-\frac{2}{h}\Phi_0}\,dL(x) \leq \BigO(h^{-\rho})||(\frakp^w(x,hD_x) - z)u||\:||u||+ \BigO(h^{1-\rho})||u||^2.
\end{equation}
Here, norms are taken in $L^2_{\Phi_0}(\Bbb{C}^d;h).$

\end{proposition}

\begin{proof}

The integral on the left-hand side of (\ref{ExteriorEstimateEqn}) is a natural candidate for the quantization-multiplication formula in Lemma \ref{QuantMult}, using
$$
	\langle \chi(h^{-\rho/2}\cdot) \frakp^w(x,hD_x) u,u\rangle_{L^2_{\Phi_0}(\Bbb{C}^d;h)}.
$$
To apply the formula, we rescale as in Appendix \ref{ChangeVarsSubsection}, using $\alpha = h^\rho$ which gives the semiclassical parameter $\tilde{h} = h/\alpha = h^{1-\rho}.$ We will use $||\cdot||_h$ for $||\cdot||_{L^2_{\Phi_0}(\Bbb{C}^d;h)}$ and likewise for inner products. Rescaling the symbol using (\ref{eRescalingSymbols}), we obtain
\begin{eqnarray*}\langle \chi(\alpha^{-1/2}x)(\frakp^w(x,hD_x)-z) u,u\rangle_h &=& \langle \chi(x)
\frakU_\alpha (\frakp^w(x,hD_x)-z) u, \frakU_\alpha u\rangle_{\tilde{h}}
\\ &=& \langle \chi \frakU_\alpha (\frakp^w(x,hD_x)-z) \frakU_\alpha^{-1}\frakU_\alpha u, \frakU_\alpha u\rangle_{\tilde{h}} \\ &=& \langle \chi
(\frakp_\alpha^w(x,\tilde{h}D_x)-z)\frakU_\alpha u, \frakU_\alpha u\rangle_{\tilde{h}}.\end{eqnarray*}

Upon finding that the natural symbol is $$\frakp_\alpha(x,\xi) = \frakp(\alpha^{1/2}x,\alpha^{1/2}\xi),$$ we then multiply by
$\alpha^{-1}$ in order to get the largest possible symbol whose second derivatives are still $\BigO(1)$ uniformly in $h.$  This follows
from the chain rule, which gives an extra power of $\alpha^{1/2}$ for every derivative and hence gives a nonnegative power of
$0<\alpha \ll 1$ for two derivatives or more.

Furthermore, this is the natural scale on which the symbol $\alpha^{-1}\frakp_\alpha(x,\xi(x))$ is uniformly bounded from below on
$\operatorname{supp}\chi.$ Recall that $\xi(x) = -\jvIm x$ so that $(x,\xi(x)) \in \Lambda_{\Phi_0},$ the natural domain of definition of
$\frakp.$ Since
$$
	\jvRe \frakp(x,\xi(x)) \geq \frac{1}{C}|(x,\xi(x))|^2 \geq \frac{1}{C}|x|^2
$$
near $(0,0)$ and is bounded from below elsewhere, for $\alpha \ll 1$ we have $\alpha^{-1}\jvRe
\frakp_\alpha(x,\xi(x)) \geq \frac{1}{C}$ uniformly in $\alpha$ on $\operatorname{supp}(\chi)\not\ni 0.$  This gives
\begin{equation}\label{echipBddBelow}
	\chi(x)\jvRe\left(\alpha^{-1}(\frakp_\alpha(x,\xi(x))-z)\right) \geq \frac{1}{C}\chi(x),
\end{equation}
as $\alpha = h^{\rho}$ for $0 < \rho < 1$ and $z$ is only logarithmically larger than $h,$ and so $\alpha^{-1}z \rightarrow 0$ as $h\rightarrow 0^+.$

Multiplying $$\langle \chi(\alpha^{-1/2}x)(\frakp^w(x,hD_x)-z) u,u\rangle_h = \langle \chi (\frakp_\alpha^w(x,\tilde{h}D_x)-z)\frakU_\alpha u,
\frakU_\alpha u\rangle_{\tilde{h}}$$ through by $\alpha^{-1}$ and applying the quantization-multiplication formula in Lemma \ref{QuantMult} gives
\begin{multline*}
	\alpha^{-1}\langle\chi(\alpha^{-1/2}x)(\frakp^w(x,hD_x)-z)u,u\rangle_h
	\\= \int_{\Bbb{C}^d}\chi(x)\alpha^{-1}(\frakp_\alpha(x,\xi(x))-z)|(\frakU_\alpha u)(x)|^2 e^{-2\Phi_0(x)/\tilde{h}}\,dL(x)+
\BigO(\tilde{h})||\frakU_\alpha u||^2_{\tilde{h}}.
\end{multline*}

We take real parts of both sides and use (\ref{echipBddBelow}) on the right-hand side.  The unitary change of variables
$\frakU_\alpha^{-1}$ undoes $\frakU_\alpha$ in the integral on the right-hand side and returns $\tilde{h}$ to $h$, and therefore
\begin{multline}\label{eExtAlmostDone}
	\alpha^{-1}\jvRe \langle\chi(\alpha^{-1/2}x)(\frakp^w(x,hD_x)-z)u,u\rangle_h
	\\ \geq \frac{1}{C} \int_{\Bbb{C}^d} \chi(\alpha^{-1/2}x)|u(x)|^2 e^{-\frac{2}{h}\Phi_0(x)}\,dL(x) - \BigO(\tilde{h})||u||^2_h.
\end{multline}
Applying the Cauchy-Schwarz inequality and discarding the cutoff function $\chi$ gives
$$\alpha^{-1}||(\frakp^w(x,hD_x)-z)u||\:||u|| \geq \alpha^{-1}\jvRe
\langle\chi(\alpha^{-1/2}x)(\frakp^w(x,hD_x)-z)u,u\rangle$$ in $L^2_{\Phi_0}(\Bbb{C}^d;h).$
Moving the $\BigO(\tilde{h})||u||^2$ to the other side of (\ref{eExtAlmostDone}) and
recalling that $\alpha = h^{\rho}$ while $\tilde{h} = h/\alpha = h^{1-\rho}$
gives the conclusion of the proposition,
(\ref{ExteriorEstimateEqn}).

\end{proof}

\section{Proof of Theorem}\label{SectionGluingProof}

We will now prove Theorem \ref{MainTheorem}. We note that Proposition \ref{UnconjugatedQuadratic}, in effect, proves Theorem \ref{MainTheorem} for the quadratic part of $p$ near $(x,\xi) = (0,0).$  To localize to the part where $p$ is nearly quadratic, we conjugated with the FBI transform and obtained the localized estimate in Proposition \ref{RescaledCutoff}, which gives an estimate like that in Theorem \ref{MainTheorem} except localized to a neighborhood of size approximately $h^{\rho/2}$ for $0 < \rho < 1.$  The price of that localization is an error living on a compact neighborhood of $\opnm{supp}(\nabla\chi_0(h^{-\rho/2}\cdot)).$  We will therefore combine the good localized estimate for $\frakq$ given by Proposition \ref{RescaledCutoff} with local approximation of $\frakq^w$ by $\frakp^w$ given by Proposition \ref{LocalFullSymbol}.  Both the error in Proposition \ref{RescaledCutoff} and the nonlocalized part of the function $u$ will be bounded using ellipticity of $\frakp$ off a neighborhood of the origin, as made precise in Proposition \ref{ExteriorEstimate}.

The computation here follows the same outline as in \cite{HeSjSt05} and in \cite{H-PS}.  Because of differing $h$-dependent scaling for the quadratic estimate in Proposition \ref{UnconjugatedQuadratic} and for the localization in Proposition \ref{RescaledCutoff}, Lemma \ref{LocalFullSymbol}, and Proposition \ref{ExteriorEstimate}, matching different powers of $h$ becomes somewhat technical.

To simplify matters, we will use some truncated notation. For cutoff functions $\chi_j(x), j = 0,1,$ introduced below,
we define rescaled versions,
$$
	\chi_j^*(x) = \chi_j(h^{-\rho/2}x),
$$
and likewise for the indicator function $\un_K.$ We write $\frakp^w$ and $\frakq^w$ for $\frakp^w(x,hD_x)$ and $\frakq^w(x,hD_x),$ where $\frakp = p \circ \kappa_T^{-1}$ and $\frakq = q \circ \kappa_T^{-1}$ are the FBI versions of the symbol $p$ and its quadratic part $q,$ satisfying all hypotheses of Theorem \ref{MainTheorem}.  All norms, unless stated otherwise, will be in $L^2_{\Phi_0}(\Bbb{C}^d;h).$  We also omit the reminder that all these statements hold for $h \in (0, h_0]$ for some $h_0 > 0$ sufficiently small.

We assume that $\gamma > 0$ and $0 < \rho < 1.$ We take a cutoff function $\chi_0\in C^\infty_0(\Bbb{C}^d)$ taking values in $[0,1]$ with $\chi_0(x) \equiv 1$ on an open set
containing $0 \in \Bbb{C}^d$.  When choosing $K$ a compact neighborhood of $\opnm{supp}(\nabla\chi_0),$ we may therefore assume that
$0 \notin K.$  Let $u \in H_{\Phi_0}(\Bbb{C}^d;h)$ and let $z \in \Bbb{C}$ obey $|z| \leq hf(h)$ and
$$
\operatorname{dist}(\operatorname{Spec} (q^w(x,hD_x)), z) \geq hf(h)^{(1-d)/2},
$$
with $f(h)$ defined as in (\ref{Definingf}) for $C$ sufficiently large to satisfy Proposition \ref{RescaledCutoff}.

Without further ado, we may begin computations.  From Proposition \ref{RescaledCutoff} and the triangle inequality, we have
\begin{multline}\label{eProofl1}
	||u|| \leq ||\chi_0^* u|| + ||(1-\chi_0^*) u||
	\\ \leq \BigO(h^{-1-\gamma(1-\rho)}) ||\chi_0^*(\frakq^w - z)u|| + \BigO(h^{-\gamma(1-\rho)})||\un^*_K u|| + ||(1-\chi_0^*)u||.
\end{multline}
Using Lemma \ref{LocalFullSymbol}, also along with the triangle inequality, we have
\begin{equation}\label{eProofl2}
	||\chi_0^*(\frakq^w - z)u|| \leq ||\chi_0^*(\frakp^w - z)u|| + ||\chi_0^*(\frakp^w - \frakq^w)u|| \leq
||\chi_0^*(\frakp^w - z)u|| + \BigO(h^{3\rho/2})||u||.
\end{equation}

To bound $||(1-\chi_0^*)u||$ and $||\un_K^*u||,$ which both are localized off a small neighborhood of zero, we introduce a cutoff function localizing near infinity.
Let $\chi_1 \in C^\infty(\Bbb{C}^d)$ take values in $[0,1]$ and let $\chi_1(x) \equiv 1$ off a small neighborhood of zero, so small that
$\chi_1(x) \equiv 1$ whenever $x \in K$ or $\chi_0(x) \neq 1.$  We furthermore
choose $\chi_1(x)$ so that, in an even smaller neighborhood of zero, $\chi_1(x) \equiv 0.$  Since $\un_K$ and $(1-\chi_0)$ are dominated by $\chi_1$ and since $h^{-\gamma(1-\rho)} \gg 1,$ we combine (\ref{eProofl1}) and (\ref{eProofl2}) to obtain
\begin{equation}\label{eProofl3}
	||u|| \leq \BigO(h^{-1-\gamma(1-\rho)})||\chi_0^*(\frakp^w - z)u|| + \BigO(h^{-1-\gamma(1-\rho) + 3\rho/2})||u|| +
\BigO(h^{-\gamma(1-\rho)})||\chi_1^* u||.
\end{equation}

The portion of $u$ localized to a neighborhood of infinity by $\chi_1^*$ can be estimated by Proposition \ref{ExteriorEstimate}.  We see that
$$
	||\chi_1^* u||^2 = \int \chi_1(h^{-\rho/2})^2 |u|^2 e^{-\frac{2}{h}\Phi_0(x)}\,dL(x) \leq \BigO(h^{-\rho})||(\frakp^w-z)u||\:||u|| + \BigO(h^{1-\rho})||u||^2.
$$
We may take square roots of both sides and distribute the square root over addition at the price of a constant.  The term with $||(\frakp^w-z)u||\:||u||$
may be expanded with the Peter-Paul form of the Cauchy-Schwarz inequality.  Therefore, for constants $\alpha \in \Bbb{R}$ and $M > 0$ to be chosen, we have
\begin{equation}\label{eProofExt}
	||\chi_1^* u|| \leq \BigO_M (h^{-\rho-\alpha})||(\frakp^w-z)u|| + \frac{1}{M}h^\alpha ||u|| + \BigO(h^{(1-\rho)/2})||u||.
\end{equation}

We insert (\ref{eProofExt}) into (\ref{eProofl3}) and rearrange terms so that all terms involving $||u||$ are on the left and all terms involving $||(\frakp^w-z)u||$ are on the right.  When rearranging, it only increases the larger side of the inequality to discard $\chi_0^*$ from $||\chi_0^*(\frakp^w-z)u||.$  We obtain
\begin{multline}\label{eProofl4}
	||u|| - \BigO(h^{-1-\gamma(1-\rho) + 3\rho/2})||u|| - \BigO(h^{-\gamma(1-\rho)})\left[\frac{1}{M}h^\alpha + \BigO(h^{(1-\rho)/2})\right]||u||
	\\ \leq \BigO(h^{-1-\gamma(1-\rho)})||(\frakp^w - z)u|| + \BigO_M (h^{-\gamma(1-\rho)-\rho-\alpha})||(\frakp^w-z)u||.
\end{multline}
We aim to absorb every term in the left-hand side into $||u||$ with a loss of at most a constant.  We therefore make the restrictions that
\begin{equation}\label{eProofRestrict1}
	-1-\gamma(1-\rho)+3\rho/2 > 0,
\end{equation}
\begin{equation}\label{eProofRestrict2}
	-\gamma(1-\rho)+(1-\rho)/2 > 0.
\end{equation}
We can also choose $M$ after having chosen $\gamma$ and $\rho,$ which fixes the $h$-independent constants from Proposition \ref{RescaledCutoff},
Lemma \ref{LocalFullSymbol}, and Proposition \ref{ExteriorEstimate}.  It is therefore sufficient to have
\begin{equation}\label{eProofRestrict3}
	-\gamma(1-\rho) + \alpha = 0,
\end{equation}
and at the end we can choose $M$ large enough but fixed.

From (\ref{eProofRestrict2}) and the fact that $\rho < 1$ we see that we must have
\begin{equation}\label{eGammaHalf}
	\gamma < \frac{1}{2}.
\end{equation}
This does obstruct any attempt to obtain, by increasing $\gamma$, a smaller implied constant $C$ for $f(h)$ in (\ref{Definingf}),
but the author is not aware of a context in which such an increase in $\gamma$ would be useful.  Instead, we are concerned with $\gamma$ small, so that the restriction (\ref{eGammaHalf})
is not a problem.  To establish (\ref{eProofRestrict1}), it suffices to choose $\rho$ such that
\begin{equation}\label{eRho34}
	1> \rho > \frac{\gamma+1}{\gamma+3/2},
\end{equation}
and from this
$$
\gamma(1-\rho) < \gamma \left(\frac{1/2}{\gamma+3/2}\right) = \frac{1}{2} - \frac{3}{4\gamma+6}.
$$
This function is increasing for $\gamma>-3/2$, and we obtain that, for any $\gamma\in (0, 1/2)$,
$$
\gamma(1-\rho) < \frac{1}{2}  - \frac{3}{8}=\frac{1}{8}.
$$

What remains is to analyze the second term on the right-hand side of (\ref{eProofl4}).
Using (\ref{eProofRestrict3}) and (\ref{eGammaHalf}), we see that
$$
	-\gamma(1-\rho)-\rho-\alpha = -2\gamma(1-\rho) - \rho > -(1-\rho)-\rho = -1.
$$
Therefore
$$
	h^{-1-\gamma(1-\rho)} \gg h^{-\gamma(1-\rho)-\rho-\alpha},
$$
and so the second term on the right-hand side of (\ref{eProofl4}) may be absorbed into the first for $h$ sufficiently small.

Choosing $\rho$ according to (\ref{eRho34}), we now have the FBI transform side estimate
$$
	||u||_{H_{\Phi_0}(\Bbb{C}^d;h)} \leq \BigO(h^{-1-\gamma(1-\rho)})||(\frakp^w(x,hD_x)-z)u||_{H_{\Phi_0}(\Bbb{C}^d;h)},
$$
for any $u \in H_{\Phi_0}(\Bbb{C}^d;h)$ and for all $h$ sufficiently small.  Conjugating with the FBI transform, which is an
isometry, we obtain the $L^2(\Bbb{R}^d)$ estimate,
\begin{equation}\label{eProofConc}
	||u||_{L^2(\Bbb{R}^d)} \leq \BigO(h^{-1-\gamma(1-\rho)})||(p^w(x,hD_x)-z)u||_{L^2(\Bbb{R}^d)}.
\end{equation}
This implies that $p^w(x,hD_x) - z$ is injective.  The fact that the resolvent therefore exists is a consequence of the fact that, in tiny neighborhoods of $0 \in \Bbb{C},$ we know that the family
$\{p^w(x,hD_x) - z\}$ is a holomorphic family of Fredholm operators of index zero.  While this is a standard fact, it has also been proven along the way in Proposition \ref{ResolventUpperBounds}, via the decomposition (\ref{eResolvents1}).

Equation (\ref{eProofConc}) therefore proves the resolvent bound (\ref{EqnMainResolventBound}), completing the proof of Theorem \ref{MainTheorem}.

\appendix

\section{Changes of Variables and the FBI-Bargmann
Transform}\label{ChangeVariablesAppendix}

In the paper, it is frequently useful to change scales both on the
real $L^2$ side and on $H_{\Phi_0},$ its image under the standard
FBI-Bargmann transform.  The formulas are valid for any
FBI-Bargmann transform with holomorphic quadratic phase, and so we begin by
recalling the definitions and properties of such transforms,
drawing heavily from lectures by J. Sj\"{o}strand
\cite{SjLoR}. We refer the interested reader to Chapters 12.2 and
12.3 there for further explanation and proofs.

After the introduction to the FBI-Bargmann transform and pseudodifferential operators on the FBI transform side appears a list of
easily-verified but convenient change of variables identities.  Next, we record a general procedure for applying Schur's test to some integral operators on the FBI transform side.
Finally, there is a proof of a quantization vs. multiplication
result which provides a way of bounding operators from below on
sets where the symbol is bounded from below.

\subsection{The FBI-Bargmann Transform}\label{aFBIIntro}

For $u \in \mathcal{S}'(\Bbb{R}^d)$ and $x \in \Bbb{C}^d$ define, for a particular
$h$-independent constant $C_T$ explained below,
\begin{equation}\label{eFBIDef}
	T_h u(x) := C_Th^{-3d/4}\int_{\Bbb{R}^d} e^{\frac{i}{h}\varphi(x,y)}u(y)\,dy,\quad 0< h \leq 1,
\end{equation}
where $\varphi(x,y)$ is a holomorphic quadratic form on $\Bbb{C}^{2d}.$  In
order to have decay in the exponential factor we insist that
$$
	\jvIm \varphi_{yy}'' = \jvIm \left(\frac{\partial^2\varphi} {\partial y_j \partial y_k} \right)_{j,k=1}^d > 0
$$
in the sense of being a positive definite matrix, and we make the further assumption that
\begin{equation}\label{CrossPartialsInvertible}
	\operatorname{det}\varphi_{xy}'' = \operatorname{det}\left(\frac{\partial^2\varphi} {\partial x_j \partial y_k} \right)_{j,k=1}^d \neq 0.
\end{equation}

Maximizing the real part of the exponent in (\ref{eFBIDef}) leads naturally to the function $$\Phi(x) = \sup_{y\in\Bbb{R}^d} -\jvIm \varphi(x,y),$$
which must be a real strictly plurisubharmonic quadratic form on $\Bbb{C}^d.$ For a suitable choice of $C_T > 0,$ it is known that $T_h$ is a unitary map from $L^2(\Bbb{R}^d)$ onto
$$H_\Phi(\Bbb{C}^d;h) := \mathrm{Hol}(\Bbb{C}^d) \cap
L^2(\Bbb{C}^d; e^{-\frac{2}{h}\Phi(x)}\,dL(x)),$$ where
$\mathrm{Hol}(\Bbb{C}^d)$ is the set of holomorphic functions on
$\Bbb{C}^d$ and $dL(x)$ is Lebesgue measure on $\Bbb{C}_x^d$ when
viewed as $\Bbb{R}^d_{\jvRe x} \times \Bbb{R}^d_{\jvIm x}.$  In
other words, $H_\Phi(\Bbb{C}^d;h)$ consists of those holomorphic
functions $u:\Bbb{C}^d \rightarrow \Bbb{C}$ for which the norm
$$||u||_{H_\Phi(\Bbb{C};h)}^2 := \int_{\Bbb{C}^d} |u(x)|^2
e^{-\frac{2}{h}\Phi(x)}\,dL(x)$$ is finite.

When modifying functions on the FBI-Bargmann side, holomorphy may not be preserved, so we also define $$L^2_\Phi(\Bbb{C}^d;h) = L^2(\Bbb{C}^d;
e^{-\frac{2}{h}\Phi(x)}\,dL(x)).$$

Associated with $\varphi$ and $\Phi$ is the submanifold of
$\Bbb{C}^{2d}$
$$\Lambda_\Phi := \left\{\left(x,\frac{2}{i}\frac{\partial \Phi}{\partial
x}\right)\::\: x\in \Bbb{C}^d\right\}$$ and the canonical transformation
$$\kappa_T: \Bbb{C}^{2d} \ni (y, -\varphi'_y(x,y))\mapsto
(x,\varphi'_x(x,y))\in \Bbb{C}^{2d}.$$  In view of our assumption
(\ref{CrossPartialsInvertible}), the map $\kappa_T$ is
well-defined.  As we are assuming $\varphi$ is holomorphic and quadratic, it
follows that $\kappa_T$ must be complex linear. The canonical
transformation is an isomorphism between real symplectic manifolds $\kappa_T: \Bbb{R}^{2d}\rightarrow \Lambda_\Phi.$

The classical Bargmann transform, which is all that is necessary
for the results in the body of the paper, is given by
$$\varphi_0(x,y) = \frac{i}{2}(x-y)^2.$$ A computation shows that
$$\Phi_0(x) = \frac{1}{2}(\jvIm x)^2,$$ and $$\kappa_T(y,\eta)
= (y-i\eta, \eta).$$

\subsection{Pseudodifferential Operators on the FBI-Bargmann Side}

There are many nearly-equivalent ways of defining the Weyl
quantization of symbols in class $$S(\Lambda_\Phi, 1) :=
\{\mathfrak{a}\in C^\infty(\Lambda_\Phi,
\Bbb{C})\::\:|\partial^\alpha \mathfrak{a}| = \BigO_\alpha(1)\},$$
where derivatives need not be holomorphic.  More details can be
found in \cite{SjLoR} (Chapter 12), \cite{HeSjSt05} (Section 3),
or \cite{H-PS} (Section 3), in descending order of detail.

We can write
$$
	\mathfrak{a}^w(x,hD_x)u(x) = \frac{1}{(2\pi h)^d} \iint_{\Gamma_0(x)}e^{\frac{i}{h}(x-y)\cdot\theta} \mathfrak{a}(\frac{x+y}{2},\theta)u(y) \,dy\,d\theta,
$$
with contour
$$
	\Gamma_0(x) := \left\{\theta = \frac{2}{i}\frac{\partial \Phi}{\partial x}\left(\frac{x+y}{2}\right)\right\}.
$$
Note that this is equivalent to taking $(\frac{x+y}{2},\theta) \in \Lambda_\Phi.$
So defined, we have a uniformly bounded map $$\mathfrak{a}^w(x,hD_x) = \BigO(1):H_\Phi\rightarrow
H_\Phi,$$ for any $\mathfrak{a}\in S(\Lambda_\Phi, 1).$ Furthermore, we have an exact version of the Egorov theorem relating the Weyl
quantization of $S(1)$ on $L^2(\Bbb{R}^d)$ and of $S(\Lambda_\Phi, 1)$ on $H_\Phi,$ stating that 
\begin{equation}\label{eExactEgorov}
	Ta^w(x,hD_x)T^{-1} = (a\circ\kappa_T^{-1})^w(x,hD_x).
\end{equation}

It is possible to shift the contour $\Gamma_0(x)$ to make the
integral converge absolutely, and to cut the integrand off away from the
diagonal $\{x = y\}$ with small error: letting $\psi_0(x)\in
C_0^\infty(\Bbb{C}^d)$ be a standard cutoff function taking values in $[0,1]$ and equal to 1 near $x = 0$,
we have $$\mathfrak{a}^w(x,hD_x)u(x) = \frac{1}{(2\pi h)^d}
\iint_{\Gamma_0(x)} e^{\frac{i}{h}(x-y)\cdot\theta}
\psi_0(x-y)\mathfrak{a}(\frac{x+y}{2},\theta)u(y)\,dy\,d\theta +
R_0u(x).$$

One can construct an almost holomorphic extension of $\mathfrak{a}$ (which
we will continue to refer to using $\mathfrak{a}$) to a tubular
neighborhood of $\Lambda_\Phi \subseteq \Bbb{C}^d$ obeying
$$|\bar{\partial}\mathfrak{a}(x,\xi)| = \BigO_N(1)|\xi - \frac{2}{i}\frac{\partial \Phi}{\partial x}(x)|^N,~~\forall N \geq 0.$$
(See Chapter 8 of \cite{DiSj} for a description of methods of construction, which may be easily adapted to our situation, and for references.)
Defining
\begin{equation}\label{eDefGammat}
	\Gamma_t(x) := \left\{\theta = \frac{2}{i}\frac{\partial \Phi}{\partial x}\left(\frac{x+y}{2}\right) + it\overline{(x-y)}\right\}
\end{equation}
for $t \geq 0$, we shift the contour $\Gamma_0$ and arrive at
\begin{multline}\label{ContourRealizationOfWeyl}\mathfrak{a}^w(x,hD_x)u(x)
= \frac{1}{(2\pi h)^d} \iint_{\Gamma_1(x)}
e^{\frac{i}{h}(x-y)\cdot\theta}
\psi_0(x-y)\mathfrak{a}(\frac{x+y}{2},\theta)u(y)\,dy\,d\theta \\
+ R_0u(x) + R_1u(x).\end{multline}  The remainder $R_1u(x)$ is computed via Stokes' theorem on $\bigcup_{0\leq t \leq 1} \Gamma_t(x).$ In
\cite{SjLoR} and \cite{HeSjSt05} it is shown that the errors are negligible up to arbitrarily high powers of $h$:
\begin{equation}\label{eWeylShiftErrors}
	R_0, R_1 = \BigO(h^\infty):L^2_{\Phi}(\Bbb{C}^d;h) \rightarrow L^2_{\Phi}(\Bbb{C}^d;h).
\end{equation}

\subsection{Changes of Variables}\label{ChangeVarsSubsection}

We use
\begin{equation}\label{eRealU}
	U_\alpha(u)(y) = \alpha^{d/4}u(\alpha^{1/2}y)
\end{equation}
and
\begin{equation}\label{eFBIU}
\frakU_\alpha(v)(x) = \alpha^{d/2}v(\alpha^{1/2}x),
\end{equation}
the natural FBI-Bargmann side
analogue.  There is a natural scaled semiclassical parameter corresponding to $U_\alpha$ or $\frakU_\alpha$ which we
will use throughout the paper: $$\tilde{h} = \frac{h}{\alpha}.$$ The constant factors in definitions (\ref{eRealU}), (\ref{eFBIU}) are
chosen so that
$$U_\alpha:L^2(\Bbb{R}^d) \rightarrow L^2(\Bbb{R}^d)$$ and
$$\frakU_\alpha: L^2_\Phi(\Bbb{C}^d; h)\rightarrow
L^2_\Phi(\Bbb{C}^d;\tilde{h})$$ are both unitary maps, and we note that $H_\Phi(\Bbb{C}^d;h)$ is mapped to $H_\Phi(\Bbb{C}^d;\tilde{h})$ by the latter.

We note for future reference that
\begin{equation}\label{ChangeVariablesDistributeProducts}U_\alpha(fg) = f(\alpha^{1/2}\cdot)U_\alpha
g,\end{equation}
and that the analogue holds for $\frakU_\alpha.$

For differential operators with quadratic symbol $q(x,\xi)$, we have
\begin{equation}\label{QuadraticChangeVariables}U_\alpha
q^w(y,hD_y)U_\alpha^{-1} = \alpha q^w(y,\tilde{h}D_y),\end{equation}
and the analogue
$$
	\frakU_\alpha \frakq^w(x,hD_x)\frakU_\alpha^{-1} = \alpha \frakq^w(x,\tilde{h}D_x)
$$
holds on the FBI transform side.

The change of variables pairs nicely with the FBI-Bargmann
transform: \begin{equation}\label{FBIChangeVariables}\frakU_\alpha T_h U_\alpha^{-1} = T_{\tilde{h}}.\end{equation}

A simple computation reveals that $$U_\alpha
p^w(y,hD_y)U_\alpha^{-1} = p_\alpha^w(y,\tilde{h}D_y),$$ where
$p_\alpha(y,\eta) = p(\alpha^{1/2}y,\alpha^{1/2}\eta).$
Conjugation with the FBI transform or direct computation yields
\begin{equation}\label{eRescalingSymbols}
	\frakU_\alpha \frakp^w(x,hD_x)\frakU_\alpha^{-1} = \frakp_\alpha^w(x,\tilde{h}D_x)
\end{equation}
with the similar $\frakp_\alpha(x,\xi) = \frakp(\alpha^{1/2}x,\alpha^{1/2}\xi).$

\subsection{Schur's Test and $H_\Phi.$}\label{SchurAppendix}

We will occasionally need to use Schur's test when confronted with an integral operator on $H_\Phi(\Bbb{C}^d;h).$  Note that
$$
	Mu(x) := e^{-\Phi(x)/h}u(x)
$$
is an isometry from $H_\Phi(\Bbb{C}^d;h)$ to a closed subspace of $L^2(\Bbb{C}^d, dL)$,
with $dL$ representing Lebesgue measure.  We may also regard $M$ as an isometric isomorphism when defined on all of $L^2_\Phi(\Bbb{C}^d;h)$.  When given an operator $A$ acting on $H_\Phi(\Bbb{C}^d;h)$ or on $L^2_\Phi(\Bbb{C}^d;h)$, it is natural to work in the more familiar space, $L^2,$ by considering $M A M^{-1}.$

If $A$ is given by an integral kernel,
$$
	Au(x) = \int_{\Bbb{C}^d} K(x,y;h)u(y)\,dL(y),
$$
it is easy to see that $M A M^{-1}$ is given by an integral kernel,
\begin{equation}\label{eConjugatedKernel}
	\tilde{K}(x,y;h) = e^{-\Phi(x)/h}K(x,y;h) e^{\Phi(y)/h}.
\end{equation}
Schur's test gives the upper bound
\begin{equation}\label{eSchurTestHPhi}
	||A||_{L^2_\Phi \rightarrow L^2_\Phi}^2 \leq \left(\sup_{x\in\Bbb{C}^d} ||\tilde{K}(x,\cdot;h)||_{L^1(\Bbb{C}^d,dL)}\right)
\left(\sup_{y\in\Bbb{C}^d} ||\tilde{K}(\cdot,y;h)||_{L^1(\Bbb{C}^d,dL)}\right).
\end{equation}

We will be most interested in applying (\ref{eSchurTestHPhi}) to integral kernels derived from the quantization formula (\ref{ContourRealizationOfWeyl}).
In this case, we will have an operator of the form
\begin{equation}\label{eGeneralWeylOperator}
	Bu(x) = (2\pi h)^{-d} \iint_{\Gamma_1} e^{\frac{i}{h}(x-y)\cdot\theta} b(x,y,\theta)u(y)\,dy\,d\theta.
\end{equation}
Recall that $\Gamma_1$ is the contour where
\begin{equation}\label{eDefThetaSchurSub}
	\theta = \frac{2}{i}(\partial_x \Phi)(\frac{x+y}{2}) + i\overline{(x-y)}.
\end{equation}

Our first step is to show that $dy\,d\theta$ is actually Lebesgue measure $dL(y)$ up to a non-vanishing constant.
Since $\Phi$ is quadratic, we note that $\partial_x\Phi$ is linear.  In fact, using symmetric inner products,
$$
	\Phi(x) = \frac{1}{2}\langle x, (\partial_x^2\Phi)x\rangle + \langle x, (\partial_x\bar{\partial}_x \Phi) \bar{x}\rangle + \frac{1}{2}\langle \bar{x}, (\bar{\partial}_x^2\Phi)\bar{x}\rangle,
$$
using that the second derivatives of $\Phi$ are constant.  We may therefore check that
$$
	\theta = (-i)\left((\partial_x\bar{\partial}_x \Phi) + I\right)\bar{y} + g(x, \bar{x}, y)
$$
for a linear function $g$.  Since we are holding $x$ constant when we integrate in $\theta,$
$$
	d\theta = \left|\det(\partial_x\bar{\partial}_x \Phi + I)\right|\,d\bar{y} + |\det(\partial_y g(x,\bar{x},y))|\,dy,
$$
with the coefficient of $dy$ immaterial because it vanishes when we integrate against $dy\,d\theta$.
Since $\Phi(x)$ is a real-valued strictly plurisubharmonic quadratic form, $\partial_x\bar{\partial}_x \Phi$ is a constant positive definite
Hermetian matrix, a fact which is unchanged upon adding $I$.  Therefore $\det(\partial_x\bar{\partial}_x \Phi + I)$ is a positive constant.
Using
$$
	dy_j = d(\jvRe y_j) + i d(\jvIm y_j), \qquad d\bar{y}_j = d(\jvRe y_j) - i d(\jvIm y_j),
$$
one can verify that
$$
	|dy \wedge d\bar{y}| = 2^d dL(y).
$$
We conclude that
$$
	dy\,d\theta = C_\Phi\, dL(y)
$$
for
$$
	C_\Phi = 2^d \left|\det(\partial_x\bar{\partial}_x \Phi + I)\right|.
$$
Consider, as an example, the classical Bargmann transform outlined at the end of Appendix \ref{aFBIIntro}.  Recalling that 
$$
	\Phi_0(x) = \frac{1}{2}(\jvIm x)^2 = -\frac{1}{8}(x-\bar{x})^2,
$$
we have that $\partial_x\bar{\partial}_x \Phi_0 = (1/4)I$ and so
$$
	C_{\Phi_0} = 2^d \opnm{det}\left(\frac{5}{4}I\right) = \left(\frac{5}{2}\right)^d.
$$

Returning to the integral operator (\ref{eGeneralWeylOperator}), we obtain a modified kernel as in (\ref{eConjugatedKernel}) given by the formula
\begin{equation}\label{eSchurWeylKernel}
	\tilde{K}(x,y) = (2\pi h)^{-d} C_\Phi\exp\left(-\frac{1}{h}(\Phi(x) - i(x-y)\cdot\theta - \Phi(y))\right)b(x,y,\theta).
\end{equation}
To analyze the absolute value of the exponential, we use the definition of $\theta$ in (\ref{eDefThetaSchurSub}) and look at the real part
of the exponent:
\begin{multline}\label{eSchurWeylPhase}
	\frac{1}{h}(-\Phi(x) + \jvRe(i(x-y)\cdot\theta) + \Phi(y))
	\\ = \frac{1}{h}\left(-|x-y|^2 - \Phi(x) + 2\jvRe((x-y)\cdot(\partial_x\Phi)(\frac{x+y}{2})) + \Phi(y)\right).
\end{multline}
Because $\Phi$ is quadratic, we have, by an application of the exact second-order Taylor expansion,
$$
\Phi(x)-\Phi(y) = \langle{(x-y),\nabla \Phi(\frac{x+y}{2})\rangle}
$$
Plugging this into (\ref{eSchurWeylPhase}) and then into (\ref{eSchurWeylKernel}), we conclude that, in this application,
\begin{equation}\label{eSchurWeylFinalUpBd}
	|\tilde{K}(x,y)| = (2\pi h)^{-d} C_\Phi e^{-\frac{1}{h}|x-y|^2} |b(x,y,\theta)|.
\end{equation}
Since this kernel has a factor which decays rapidly as $h\rightarrow 0^+$ away from the main diagonal $\{x = y\},$ the corresponding operator is often simple to bound by Schur's test.
\subsection{Quantization vs. Multiplication}

Finally, we have a relation between quantization and
multiplication on the FBI transform side on the level of inner products
that, in the proof of the theorem, provides the critical elliptic
estimate in a region away from the doubly characteristic point at
the origin.  The proof is well-established in \cite{Sj90}, \cite{HeSjSt05}, and \cite{H-PS}, and is included here on
the grounds that it is not as well-known as some other facts about the FBI transform side.

\begin{lemma}\label{QuantMult}
Let $\frakp:\Lambda_{\Phi}\rightarrow \Bbb{C}$ be smooth and bounded with all derivatives for each fixed $h>0$. We assume that, uniformly in $h>0$,
\begin{equation}\label{eFrakpBddAtZero}
	|\frakp(0,0)| = \BigO(1),\quad |\nabla \frakp(0,0)| = \BigO(1),
\end{equation}
and
\begin{equation}\label{eFrakpHessBdd}
	|\partial^\alpha \frakp| = \BigO_\alpha(1), \quad \forall |\alpha| \geq 2.
\end{equation}
We make the further assumption that, for some fixed $N,$
\begin{equation}\label{eFrakpNotCrazy}
	|\partial^\alpha \frakp| = \BigO_\alpha(h^{-N}), \quad \forall |\alpha| \leq 1.
\end{equation}
Let $\frakp$ also denote an almost holomorphic extension off $\Lambda_\Phi$ which obeys the same estimates.  Let $\chi \in C_b^\infty(\Bbb{C}^d)$ with
\begin{equation}\label{eGradChiBdd}
	\nabla \chi(x) = \BigO(\langle x \rangle^{-1}).
\end{equation}
Define $\xi(x) = \frac{2}{i}\partial_x\Phi(x).$  Then, for any $u \in H_\Phi(\Bbb{C}^d; h)$, we have
\begin{equation}\label{eQMConclusion}
	\langle \chi\frakp^w(x, hD_x)u, u\rangle = \int \chi(x)\frakp(x,\xi(x))|u(x)|^2e^{-2\Phi(x)/h}\,dL(x) + \BigO(h)||u||^2,
\end{equation}
with inner products and norms in $L^2_\Phi(\Bbb{C}^d;h).$
\end{lemma}

\begin{proof} For future use, we reformulate the conclusion (\ref{eQMConclusion}) entirely in terms of the inner product on
$L^2_\Phi(\Bbb{C}^d;h),$ obtaining
\begin{equation}\label{eQMConclusionIP}
	\langle \chi\frakp^w(x,hD_x) u, u\rangle = \langle \chi\frakp(x,\xi(x))u, u \rangle + \BigO(h)||u||^2.
\end{equation}
In this sense, when taking inner products, the Weyl quantization is multiplication by the symbol to first order.

A brief summary of the proof is as follows: first, we Taylor expand the symbol $\frakp$ as it appears in (\ref{ContourRealizationOfWeyl})
to second order. The first order term becomes the multiplication by $\frakp(x,\xi(x))$ which is the principal term in (\ref{eQMConclusion}).
The second order terms either vanish or are shown to be bounded by $\BigO(h)||u||^2$ after integration by parts.  The error in the Taylor expansion contributes
at most $\BigO(h)||u||^2$ as a consequence of Schur's test.

It is sufficient to consider a dense set of $u \in H_\Phi(\Bbb{C}^d;h),$ as both sides of (\ref{eQMConclusion}) are continuous on $L^2_\Phi(\Bbb{C}^d;h)$ for
$u\in H_{\Phi}(\Bbb{C}^d;h).$  We therefore will assume that $\langle x \rangle^N D_x^\alpha u \in L^2_\Phi(\Bbb{C}^d;h)$ for all
$N$ and $\alpha.$  As explained in Proposition \ref{PropFBIQuadraticEquivalence} and the remark following, it is sufficient to assume that, for all $N \geq 1$, we have
$\langle x \rangle^N u \in L^2_\Phi(\Bbb{C}^d;h)$. Such $u$ are well-known to be dense in $H_\Phi(\Bbb{C}^d;h)$ as they are precisely the elements of
$T_h(\mathcal{S}(\Bbb{R}^d))$, as seen in \cite{SjLoR}, Chapter 12.

The assumption (\ref{eFrakpNotCrazy}) allows us to employ the shifted contour for $\frakp^w(x,hD_x)$ appearing in (\ref{ContourRealizationOfWeyl}), because $h^N\frakp \in S(\Lambda_\Phi,1)$ uniformly in $h.$ Therefore the errors in (\ref{eWeylShiftErrors}) are of the form $\BigO(h^{-N}h^\infty) = \BigO(h^\infty)$ and are negligible.

We Taylor expand $\frakp(\frac{x+y}{2},\theta)$ appearing in (\ref{ContourRealizationOfWeyl}) at $(x,\xi(x))\in\Lambda_\Phi.$  Because the extension off $\Lambda_\Phi$ is almost holomorphic, we have $\bar{\partial}\frakp(x,\xi(x)) = 0$ for all $x \in \Bbb{C}^d,$ so the expansion involves only holomorphic derivatives.  We thus have
\begin{equation}\label{eTaylorFrakp}
	\frakp(\frac{x+y}{2},\theta) = \frakp(x,\xi(x)) + (\partial_x\frakp)(x,\xi(x))\cdot(\frac{y-x}{2}) + (\partial_\theta\frakp)(x,\xi(x))\cdot(\theta-\xi(x)) + \mathfrak{r}(x,y,\theta).
\end{equation}
Dot products continue to represent the standard non-Hermetian dot product on $\Bbb{C}^d,$ with $v \cdot w = \sum_{j=1}^d v_j w_j.$  Because second derivatives of $\frakp$ are uniformly bounded,
\begin{equation}\label{erFirstBound}
	|\frakr(x,y,\theta)| = \BigO(|x-y|^2 + |\theta-\xi(x)|^2).
\end{equation}

Along the contour for integration in (\ref{ContourRealizationOfWeyl}), where
\begin{equation}\label{eQMTheta}
	\theta = \theta(x,y) = \frac{2}{i}(\partial_x\Phi)(\frac{x+y}{2}) + i\overline{(x-y)},
\end{equation}
we have
$$
	\theta-\xi(x) = \frac{2}{i}\left[ (\partial_x\Phi)(\frac{x+y}{2}) - (\partial_x\Phi)(x) \right] + i\overline{(x-y)}.
$$
Since $\Phi(x)$ is quadratic, $\partial_x\Phi$ is a constant (real-)linear map.  We may therefore conclude that
$$
	|\theta - \xi(x)| = \BigO(|x-y|),
$$
and this with (\ref{erFirstBound}) gives
\begin{equation}\label{erSecondBound}
	|\frakr(x,y,\theta)| = \BigO(|x-y|^2).
\end{equation}

We use (\ref{eTaylorFrakp}) to expand (\ref{ContourRealizationOfWeyl}) for $\frakp^w(x,hD_x)u(x)$ into four terms, recalling that $\psi_0$ is a smooth compactly supported function equal to 1 in a neighborhood of $0 \in \Bbb{C}^d$:
\begin{eqnarray*}
	\frakp^w(x,hD_x)u(x) & = & (2\pi h)^{-d} \frakp(x,\xi(x)) \iint_{\Gamma_1} e^{\frac{i}{h}(x-y)\cdot\theta}\psi_0(x-y) u(y)\,dy\,d\theta \\
		& & + (2\pi h)^{-d} \sum_{j=1}^d (\partial_{x_j}\frakp)(x,\xi(x))\iint_{\Gamma_1} e^{\frac{i}{h}(x-y)\cdot\theta}\psi_0(x-y)(\frac{x_j - y_j}{2}) u(y)\,dy\,d\theta \\
		& & + (2\pi h)^{-d} \sum_{j=1}^d (\partial_{\theta_j}\frakp)(x,\xi(x))\iint_{\Gamma_1} e^{\frac{i}{h}(x-y)\cdot\theta} \psi_0(x-y)(\theta_j - \xi_j(x))u(y)\,dy\,d\theta
		\\ & & + (2\pi h)^{-d} \iint_{\Gamma_1} e^{\frac{i}{h}(x-y)\cdot\theta}\psi_0(x-y)\frakr(x,y,\theta) u(y)\,dy\,d\theta.
\end{eqnarray*}
The first three lines may be dealt with using standard facts, familiar from the Fourier transform, which are established in \cite{SjLoR}, Section 12.2.  We state these facts using the language of Weyl quantizations, ignoring $\BigO(h^\infty)$ errors which arise from shifting contour and including $\psi_0$ in (\ref{ContourRealizationOfWeyl}).  First, $1^w = \operatorname{id}_{H_\Phi \rightarrow H_\Phi},$ and so the integral in the first line collapses to $u(x).$  The integrals in the second line collapse to $\frac{1}{2}(x_j - x_j)u(x) = 0$ given our assumption that that $y_j u(y) \in H_\Phi(\Bbb{C}^d;h).$  We also have that $\theta_j^w = hD_{x_j},$ with a holomorphic derivative.  We now have the simpler expansion
\begin{multline*}
	\frakp^w(x,hD_x)u(x) = \frakp(x,\xi(x))u(x) + \sum_{j=1}^d (\partial_{\theta_j}\frakp)(x,\xi(x))(hD_{x_j}-\xi_j(x))u(x)
		\\ + (2\pi h)^{-d} \iint_{\Gamma_1} e^{\frac{i}{h}(x-y)\cdot\theta}\psi_0(x-y)\frakr(x,y,\theta) u(y)\,dy\,d\theta.
\end{multline*}

We will now use Schur's test, as in Appendix \ref{SchurAppendix}, to bound the integral with $\frakr.$  We obtain an integral operator of type
(\ref{eGeneralWeylOperator}), with
$$
	b(x,y,\theta) = \psi_0(x-y)\frakr(x,y,\theta).
$$
Combining the bound (\ref{erSecondBound}), the fact that $\psi_0 \in L^\infty(\Bbb{C}^d)$, and the bound for the integral kernel $\tilde{K}(x,y)$ in (\ref{eSchurWeylFinalUpBd}), we obtain
$$
	|\tilde{K}(x,y)| \leq \BigO(1)h^{-d}|x-y|^2 e^{-|x-y|^2/h}.
$$
A change of variables shows that
$$
	||\tilde{K}(x,y)||_{L^1(\Bbb{C}^d, dL(x))} \leq \BigO(h), \quad ||\tilde{K}(x,y)||_{L^1(\Bbb{C}^d, dL(y))} \leq \BigO(h).
$$
Schur's test therefore shows that the integral operator given by $\tilde{K}$ acting on $L^2(\Bbb{C}^d, dL)$ is uniformly bounded by $\BigO(h),$ and hence the corresponding operator on $L^2_\Phi(\Bbb{C}^d;h)$ is similarly bounded.  Hence
$$
	||(2\pi h)^{-d} \iint_{\Gamma_1} e^{\frac{i}{h}(x-y)\cdot\theta}\psi_0(x-y)\frakr(x,y,\theta) u(y)\,dy\,d\theta||_{L^2_\Phi(\Bbb{C}^d;h)} \leq \BigO(h)||u||_{H_\Phi(\Bbb{C}^d;h)}.
$$

Thus far, we have that
$$
	\frakp^w(x,hD_x)u(x) = \frakp(x,\xi(x))u(x) + \sum_{j=1}^d(\partial_{\theta_j}\frakp)(x,\xi(x))(hD_{x_j} - \xi_j(x))u(x) + \BigO(h||u||_{H_\Phi(\Bbb{C}^d;h)}),
$$
where the error term is a function in $L^2_\Phi(\Bbb{C}^d;h).$  We now multiply by $\chi$, take the $L^2_\Phi(\Bbb{C}^d;h)$ inner product with $u$, and obtain
$$
	\langle \chi\frakp^w(x,hD_x)u, u\rangle = \langle \chi\frakp(x,\xi(x))u, u\rangle + \sum_{j=1}^d \langle \chi(\partial_{\theta_j}\frakp)(x,\xi(x))(hD_{x_j} - \xi_j(x))u,u\rangle + \BigO(h)||u||^2,
$$
having used the Cauchy-Schwarz inequality for the error term.  To obtain the conclusion (\ref{eQMConclusionIP}), it suffices to show that
\begin{multline}\label{eLastQMTerm}
	\langle \chi(\partial_{\theta_j}\frakp)(x,\xi(x))(hD_{x_j} - \xi_j(x))u,u\rangle
	\\ = \int_{\Bbb{C}^d} \chi(x)(\partial_{\theta_j}\frakp)(x,\xi(x))\left[(hD_{x_j} - \xi_j(x))u(x)\right]\overline{u(x)}e^{-2\Phi(x)/h} \, dL(x) = \BigO(h)||u||^2,
\end{multline}
which we will now do.

We have assumed that $ue^{-\Phi/h}$ and its derivatives decay sufficiently quickly that integration by parts, $\int (D_{x_j}v)w\,dL = -\int v(D_{x_j}w)\,dL,$ is valid.  Thus
\begin{multline}\label{ePartsOnLastQM}
	\int_{\Bbb{C}^d} \chi(x)(\partial_{\theta_j}\frakp)(x,\xi(x))\left[(hD_{x_j})u(x)\right]\overline{u(x)}e^{-2\Phi(x)/h} \, dL(x)
	\\ = -\iint_{\Bbb{C}^d} u(x) hD_{x_j}\left[\chi(x)(\partial_{\theta_j}\frakp)(x,\xi(x))\overline{u(x)}e^{-2\Phi(x)/h}\right]\,dL(x).
\end{multline}
We will treat the cases where $D_{x_j}$ hits each of these multiplicands in turn.

From (\ref{eFrakpBddAtZero}), (\ref{eFrakpHessBdd}), and the fact that $|\xi(x)| = \BigO(|x|)$ because $\Phi$ is quadratic, we have that
$$
	(\partial_{\theta_j}\frakp)(x,\xi(x)) = \BigO(\langle x \rangle).
$$
Combining this with (\ref{eGradChiBdd}), we see that
$$
	(D_{x_j}\chi(x))(\partial_{\theta_j}\frakp)(x,\xi(x)) = \BigO(1)
$$
uniformly in $h,$ giving a contribution of $\BigO(h)||u||^2$ to (\ref{ePartsOnLastQM}) when the derivative hits $\chi(x).$

Using the fact that $\chi \in C^\infty_b,$ the assumption (\ref{eFrakpHessBdd}), and the fact that that derivatives of $\xi(x)$ are second derivatives of the quadratic function $\Phi,$ we also see that
$$
	\chi(x)D_{x_j}(\partial_{\theta_j}\frakp)(x,\xi(x)) = \BigO(1)
$$
uniformly in $h.$  Therefore the contribution to (\ref{ePartsOnLastQM}) when $D_{x_j}$ hits $\partial_{\theta_j}\frakp$ is also $\BigO(h)||u||^2.$

Because $u(x)$ is holomorphic, $D_{x_j}\overline{u(x)}=0.$

Finally,
$$
	hD_{x_j}e^{-2\Phi(x)/h} = -\frac{2}{i}(\partial_{x_j}\Phi)(x)e^{-2\Phi(x)/h} = -\xi_j(x)e^{-2\Phi(x)/h}.
$$
When taken with the minus sign from integration by parts in (\ref{ePartsOnLastQM}), this exactly cancels with the $-\xi_j(x)$ in (\ref{eLastQMTerm}).

Having shown the bound $\BigO(h)||u||^2$ for (\ref{eLastQMTerm}) for all $j = 1,\dots, d,$ we have completed the proof of (\ref{eQMConclusion}).

\end{proof}

\section{Quadratic Weights and Comparison between Real and FBI
Sides}\label{WeightsAppendix}

In preliminary proofs of the main result of this work, some weighted estimates implicit in \cite{HeSjSt05} and \cite{H-PS} were used.  While those facts proved
unnecessary for the proof of the theorem, it is hopefully of independent interest to make them explicit here.

One may interpret the following propositions, \ref{PropQuadraticGraphNormEquivalence} and \ref{PropFBIQuadraticEquivalence}, in terms of the graph norm $||\cdot||_{\Gamma(A)}$ corresponding to a linear function $A:\mathcal{H}\rightarrow\mathcal{K}$ between Hilbert spaces.  When $u \in \mathcal{D}(A)$, the domain of definition for $A$, we write
$$
	||u||_{\Gamma(A)} = \sqrt{||u||_{\mathcal{H}}^2 + ||Au||_{\mathcal{K}}^2}.
$$
As usual, the square root may be distributed over the sum at the price of a constant.  Write $H$ to denote the semiclassical harmonic oscillator
$$
	H = H(x,hD_x) := x^2 + (hD_x)^2
$$
acting on $\mathcal{H} = L^2(\Bbb{R}^d)$, and let $q$ be a quadratic form which is elliptic in the sense of (\ref{eDefEllipticQuad}).  Proposition \ref{PropQuadraticGraphNormEquivalence} states that, for $u \in \mathcal{S}(\Bbb{R}^d)$, the graph norms given by $\frac{1}{h}H$ and by $\frac{1}{h}q^w(x,hD_x)$ are comparable, with $h$-independent constants:
$$
	||u||_{\Gamma(\frac{1}{h}q^w(x,hD_x))} \sim ||u||_{\Gamma(\frac{1}{h}H)}.
$$
Here, $a \sim b$ indicates that there exists an $h$-independent constant $C > 0$ with
$\frac{1}{C}a \leq b \leq Ca.$

In Proposition \ref{PropFBIQuadraticEquivalence}, the same result on the FBI transform side is proven, with multiplication by $|x|^2$ replacing the harmonic oscillator $H$:
$$
	||u||_{\Gamma(\frac{1}{h}q^w(x,hD_x))} \sim ||(1+\frac{1}{h}|x|^2)Tu||_{L^2_\Phi(\Bbb{C}^d;h)} \sim ||Tu||_{\Gamma(\frac{1}{h}|x|^2)},
$$
where the right-hand comparison is a straightforward consequence of the Cauchy-Schwarz inequality.  While this gives information on decay of $Tu \in H_\Phi(\Bbb{C}^d;h)$, we note that multiplication by $|x|^2$, which does not preserve holomorphic functions, is not an operator on $H_\Phi(\Bbb{C}^d;h)$.

One could refer instead to the (unbounded) holomorphic multiplication operator
\begin{equation}\label{ex2Sum}
	x^2 := \sum_{j=1}^d x_j^2,
\end{equation}
which does take $H_\Phi(\Bbb{C}^d;h)$ to itself while giving the same graph norm as $|x|^2$.  However, we caution that $x^2$ fails to be self-adjoint on $H_\Phi(\Bbb{C}^d;h)$.  Thus, neither FBI-side operator enjoys precisely the same pleasant properties as the harmonic oscillator $H$, but multiplication by $|x|^2$ is certainly useful and is sufficient for applications in this paper.

We now formally state and prove the comparisons described earlier, beginning with quadratic differential operators on the real side, $L^2(\Bbb{R}^d).$

\begin{proposition}\label{PropQuadraticGraphNormEquivalence}

Let $q(x,\xi)$ be an elliptic quadratic form on $\Bbb{R}^{2d}.$  Then for $u \in \mathcal{S}(\Bbb{R}^d)$ and any $h \in (0, \infty),$ we have
\begin{equation}\label{RealQuadraticWeightEquivalenceEqn}||u|| + \frac{1}{h}||q^w(x,hD_x)u|| \sim ||u|| + \frac{1}{h}||(x^2 +
(hD_x)^2)u||,\end{equation} with norms taken in $L^2(\Bbb{R}^d).$  Here $a \sim b$ means that there exists an $h$-independent constant $C > 0$ where
$\frac{1}{C}a \leq b \leq Ca.$

\end{proposition}

\begin{proof}

Norms throughout the proof will be assumed to be in $L^2(\Bbb{R}^d)$. We may first reduce to the case $h = 1.$  The change of variables
(\ref{QuadraticChangeVariables}) gives
$$
	||u|| + \frac{1}{h}||q^w(x,hD_x)u|| = ||U_h u|| + ||q^w(x,D_x) U_h u||,
$$ and similarly
$$
||u|| + \frac{1}{h}||(x^2 + (hD_x)^2)u|| = ||U_h u|| + ||(x^2 + D_x^2)U_h u||,
$$
from which the $h$-dependent equivalence (\ref{RealQuadraticWeightEquivalenceEqn}) follows immediately from the case $h = 1.$

For the $h = 1$ case, write $a^w = a^w(x,D_x)$ and consider symbols in classes $$S_{1,-1}^m(\langle (x,\xi)\rangle) := \{p: \Bbb{R}^{2d}
\rightarrow \Bbb{C} \::\: |\partial^\alpha_x\partial^\beta_\xi p(x,\xi)| = \BigO_{\alpha,\beta}(1)\langle(x,\xi)\rangle^{m-|\alpha|-|\beta|}\},$$
as usual using $\langle (x,\xi)\rangle = (1 + |x|^2 + |\xi|^2)^{1/2}.$  For simplicity, denote these classes $S^m$ for the remainder of this
proof.  Via the pseudodifferential calculus, noting that both derivatives in $x$ and $\xi$ remove one power of $\langle (x,\xi)\rangle,$ we have
that, when $a \in S^m$ and $b \in S^{m'},$ $$a^w\circ b^w - (ab)^w = c^w$$ for $c \in S^{m+m'-2}.$

For any $\chi \in C_0^\infty(\Bbb{R}^{2d})$ a standard cutoff function taking values in $[0,1]$ with $\chi \equiv 1$ near 0, note that $$\frac{1-\chi}{q} \in S^{-2}.$$ By the
pseudodifferential calculus,
$$\left(\frac{1-\chi}{q}\right)^w q^w + \chi^w - r^w = 1_{L^2 \rightarrow L^2},$$ for some $r \in S^{-2}.$  Therefore 
\begin{equation}\label{eCalculusHQ}
	||(x^2 + D_x^2)u|| \leq ||(x^2 + D_x^2)\left(\frac{1-\chi}{q}\right)^w q^w u|| + ||(x^2 + D_x^2)\chi^w u|| + ||(x^2 + D_x^2)r^w u||.
\end{equation}
Because $x^2 + \xi^2 \in S^2,$ $$(x^2 +
D_x^2)\left(\frac{1-\chi}{q}\right)^w = (\tilde{r})^w$$ for $\tilde{r} \in S^0.$  Thus, by the Calder\'{o}n-Vaillancourt theorem, $\tilde{r}^w$ is bounded from $L^2(\Bbb{R}^d)$ to $L^2(\Bbb{R}^d)$
and so the first term in (\ref{eCalculusHQ}) is bounded by $$||(x^2 + D_x^2)\left(\frac{1-\chi}{q}\right)^w q^w u|| \leq \BigO(1)||q^w u||.$$  The second and third terms
are bounded by $\BigO(1)||u||$ because $\chi \in S^{-\infty}$ and $r \in S^{-2},$ and so $(x^2 + D_x^2)\chi^w$ and $(x^2 + D_x^2)r^w$ are given by
quantizations of symbols in $S^{-\infty}$ and $S^0,$ respectively.  Therefore $$||(x^2 + D_x^2)u|| \leq \BigO(1)||u|| + \BigO(1)||q^w u||.$$

Exchanging the roles of $q$ and $x^2 + D_x^2 = (x^2 + \xi^2)^w$ above gives the reverse inequality, $$||q^w u|| \leq \BigO(1)||u|| +
\BigO(1)||(x^2 + D_x^2)u||,$$ proving the proposition.

\end{proof}

\begin{remark}\label{RemRealQuadraticDomain}

It is proven, for instance, following Corollary 4.2 of \cite{Ho95}, that the maximal realization of $q^w(x,hD_x)$ is the graph closure of the restriction of $q^w(x,hD_x)$
to $\mathcal{S}(\Bbb{R}^d).$  Therefore the estimate (\ref{RealQuadraticWeightEquivalenceEqn}) extends to $u \in \operatorname{Dom}(q^w(x,hD_x))$
which is precisely $$\operatorname{Dom}(q^w) = \{u\in L^2(\Bbb{R}^d)\::\: q^w u \in L^2(\Bbb{R}^d)\}.$$

\end{remark}

In the next proposition, we establish a similar characterization of the graph norm on $\operatorname{Dom}(q^w)$ in terms of the FBI-Bargmann transform.

\begin{proposition}\label{PropFBIQuadraticEquivalence}

Let $q(x,\xi)$ an elliptic quadratic form on $\Bbb{R}^{2d}$ and recall $T_h,$ the $h$-dependent FBI-Bargmann transform corresponding to
$\varphi(x,y)$ a valid quadratic weight (see Appendix \ref{ChangeVariablesAppendix}).  Then for any $u \in \mathcal{S}(\Bbb{R}^d),$ uniformly in
$h,$
$$||u||_{L^2(\Bbb{R}^d)} + \frac{1}{h}||q^w(x,hD_x)u||_{L^2(\Bbb{R}^d)} \sim ||\left(1+\frac{|x|^2}{h}\right)T_h u||_{L^2_\Phi(\Bbb{C}^d;h)}.$$

\end{proposition}

\begin{proof}

Again, we may reduce to the case $h = 1$ by applying (\ref{QuadraticChangeVariables}) and (\ref{FBIChangeVariables}). This gives the equivalent
statement,
\begin{equation}\label{eFBINonSCQuadEquiv}
	||U_h u||_{L^2(\Bbb{R}^d)} + ||q^w(y,D_y)U_h u||_{L^2(\Bbb{R}^d)} \sim ||(1 + |x|^2)T_1 U_h u||_{L^2_\Phi(\Bbb{C}^d;1)}.
\end{equation}
We write
$u$ instead of $U_h u$ for brevity.  Furthermore,
\begin{multline}\label{eExpJBx}
	||(1+|x|^2)T_1 u||_{L^2_\Phi(\Bbb{C}^d;1)}^2 = \int(1+|x|^2)^2|T_1 u|^2 e^{-2\Phi}\,dL(x) \\
= \int |T_1 u|^2 e^{-2\Phi}\,dL(x) + 2\int |xT_1 u|^2 e^{-2\Phi}\,dL(x) + \int |x|^4|T_1 u|^2 e^{-2\Phi}\,dL(x) \\ \sim ||T_1
u||_{H_\Phi(\Bbb{C}^d;1)}^2 + ||x^2 T_1 u||_{H_\Phi(\Bbb{C}^d;1)}^2,
\end{multline} after an application of the Cauchy-Schwarz inequality.

One inequality in the relation (\ref{eFBINonSCQuadEquiv}) follows upon showing that $$||x^2 T_1 u||_{H_\Phi(\Bbb{C}^d;1)} \leq \BigO(1)
(||u||_{L^2(\Bbb{R}^d)} + ||q^w(y,D_y)u||_{L^2(\Bbb{R}^d)}).$$  The exact Egorov theorem (\ref{eExactEgorov}), which is valid for symbol classes $S(\Lambda_\Phi, m),$ combined with the fact that the canonical transformation $\kappa_T$
is linear, allows us to express the unbounded multiplication operator $u(x)\mapsto x_k u(x)$, which acts on $H_\Phi(\Bbb{C}^d;1)$, on the $L^2(\Bbb{R}^d)$ side: $$T_1^{-1}x_kT_1 = \sum_{j=1}^d (\lambda_{jk} y_j + \mu_{jk} D_{y_j})$$ for scalars $\lambda_{jk},\mu_{jk} \in \Bbb{C}.$  Therefore $$T_1^{-1}x^2T_1 = \sum_{|\alpha| + |\beta| \leq 2}
a_{\alpha \beta} y^\alpha \partial_y^\beta$$ with $a_{\alpha\beta} \in \Bbb{C}.$ The arguments in the proof of Proposition
\ref{PropQuadraticGraphNormEquivalence} show that $$||y^\alpha \partial_y^\beta u||_{L^2(\Bbb{R}^d)} \leq \BigO(1)(||u||_{L^2(\Bbb{R}^d)} + ||q^w(y, D_y) u||_{L^2(\Bbb{R}^d)})$$ for $|\alpha| +
|\beta| \leq 2$ and any elliptic quadratic form $q.$  We may conclude from (\ref{eExpJBx}) that
$$
	||(1+|x|^2)T_1 u||_{L^2(\Bbb{C}^d;1)} \leq \BigO(1) \left(||u||_{L^2(\Bbb{R}^d)} + ||q^w(y,D_y)u||_{L^2(\Bbb{R}^d)}\right).
$$

In the other direction, the exact Egorov theorem and Fourier inversion suffice to show that $T_1 q^w(y,D_y) T_1^{-1}$ is given by a degree 2 polynomial in $\{x_j, D_{x_j}\}_{j=1}^d$.  That is,
$$
	T_1 q^w(y,D_y)T_1^{-1} = \sum_{|\alpha| + |\beta| \leq 2} b_{\alpha \beta} x^\alpha D_x^\beta,
$$
for some $b_{\alpha\beta} \in \Bbb{C}.$
Using $||[x_j,D_{x_j}]u||_{L^2_\Phi(\Bbb{C}^d;1)} = ||u||_{L^2_\Phi(\Bbb{C}^d;1)},$ it is therefore sufficient to show that
\begin{equation}\label{eMultBdsDerivs}
	||D_{x_j} u||_{L^2_\Phi(\Bbb{C}^d;1)} \leq \BigO(1) ||\langle x \rangle u||_{L^2_\Phi(\Bbb{C}^d;1)}
\end{equation}
for any $j,$ which we will now do.

Cauchy's integral formula with $\mathbf{e}_j$ the standard basis vector with $1$ in the $j$th position gives, for any $v \in H_\Phi(\Bbb{C}^d;1),$
$$|D_{x_j}v(x)|^2 \leq \BigO(1) \left(\int_0^{2\pi} \frac{|v(x+R_x e^{i\theta}\mathbf{e}_j)|}{R_x}\,d\theta\right)^2 \leq
\frac{\BigO(1)}{R_x^2}\int_0^{2\pi}|v(x+R_xe^{i\theta})|^2\,d\theta,$$ where the radius $R_x$ will be chosen depending on $x.$  The square of the
$L^2_\Phi(\Bbb{C}^d;1)$ norm is given by integrating after multiplying by $e^{-2\Phi(x)},$ and notice that, because $\Phi$ is quadratic, $$|\Phi(y) - \Phi(x)| \leq \BigO(1)
\quad\mathrm{if}\quad |y-x|\leq \BigO\left(\frac{1}{\langle x \rangle}\right).$$  Therefore, to change variables into $$w = x +
R_xe^{i\theta}\mathbf{e}_j$$ with $$e^{-2\Phi(x)} \leq \BigO(1)e^{-2\Phi(w)},$$ use $$R_x = \frac{1}{C\langle x \rangle}$$ for $C > 1$ sufficiently
large that the Jacobian determinant is bounded above and below:
$$\left|\frac{\partial{w}}{\partial{x}}\right| \sim 1.$$  The final fact we need is that $R_x \sim R_w,$ which follows from $|x - w| \leq 1.$

Changing variables from $x$ to $w$ and interchanging the order of
integration gives \begin{eqnarray*} ||D_{x_j}v||^2_{H_\Phi(\Bbb{C}^d;1)} &=&
\int_{\Bbb{C}^d} |D_{x_j} v(x)|^2 e^{-2\Phi(x)}\,dL(x) \\ &\leq&
\BigO(1)\int_0^{2\pi} \int_{\Bbb{C}^d} \frac{|v(w)|^2}{R_w^2}
e^{-2\Phi(w)}\,dL(w)\,d\theta \\ &\leq& \BigO(1)||\langle x\rangle
v||_{L^2_\Phi(\Bbb{C}^d;1)}^2,
\end{eqnarray*} on account of the fact that $R_w^{-2} \sim \langle
w \rangle^2.$

Having established that multiplication bounds derivatives in the sense of (\ref{eMultBdsDerivs}), the terms given by $||T_1u||_{H_\Phi(\Bbb{C}^d;1)}$ and $||x^2 T_1 u||_{H_\Phi(\Bbb{C}^d;1)}$ dominate all others given by the exact
Egorov theorem described previously.  This completes the proof that $$||u||_{L^2(\Bbb{R}^d)} + ||q^w(y,D_y)u||_{L^2(\Bbb{R}^d)} \leq
\BigO(1)||(1+|x|^2)T_1u||_{L^2_\Phi(\Bbb{C}^d;1)},$$ which completes the proof of the proposition.

\end{proof}

\begin{remark}

Having seen that, for elliptic $q:\Bbb{R}^{2d}\rightarrow \Bbb{C},$ we have
$$||u||_{L^2(\Bbb{R}^d)} + ||q^w u||_{L^2(\Bbb{R}^d)} < \infty$$
if and only if 
$$||(1 + |x|^2)T_h u||_{L^2_\Phi(\Bbb{C}^d;h)} < \infty,$$
the characterization 
$$\operatorname{Dom}(\frakq^w(x,hD_x)) = \{u \in H_{\Phi}(\Bbb{C}^d;h)\::\: (1+|x|^2)u \in L^2_{\Phi}(\Bbb{C}^d;h)\}$$
is proven in view of Remark \ref{RemRealQuadraticDomain}.

We also note that we may make an analogous definition of elliptic quadratic form, taken from (\ref{eDefEllipticQuad}) but adapted to the FBI transform side.  If $\frakq:\Bbb{C}^{2d}\rightarrow \Bbb{C}$ is a holomorphic quadratic form, we say that $\frakq$ is elliptic along $\Lambda_\Phi$ if $\frakq$ fulfills
\begin{equation}\label{eDefEllipticFBI}
	\frakq^{-1}(\{0\})\cap \Lambda_\Phi = \{0\}, \quad \frakq(\Lambda_\Phi) \neq \Bbb{C}.
\end{equation}
Because the linear complex canonical transformation $\kappa_T$ in the exact Egorov relation (\ref{eExactEgorov}) is a real-linear map from $\Bbb{R}^{2d}$ onto $\Lambda_\Phi$, we see that $q = \frakq\circ\kappa_T$ is an elliptic quadratic form in the sense of (\ref{eDefEllipticQuad}) if and only if $\frakq$ is an elliptic quadratic form along $\Lambda_\Phi$ in the sense of (\ref{eDefEllipticFBI}).  To see this, we begin by noting that quadratic forms on $\Bbb{R}^{2d}$ may be uniquely holomorphically extended to $\Bbb{C}^{2d}$, and complex linearity of $\kappa_T$ gives that $q$ is a holomorphic quadratic form on $\Bbb{C}^{2d}$ if and only if $\frakq$ is a holomorphic quadratic form on $\Bbb{C}^{2d}$.  Bijectivity of $\kappa_T:\Bbb{R}^{2d}\rightarrow \Lambda_\Phi$ then gives a bijection between $q^{-1}(\{0\})$, considered as a subset of $\Bbb{R}^{2d}$, and $\frakq^{-1}(\{0\}) \cap \Lambda_\Phi$. The same bijectivity ensures that $q(\Bbb{R}^{2d}) = \frakq(\Lambda_\Phi)$, completing the proof.

Finally, with the characterization $$\mathcal{S}(\Bbb{R}^d) = \{u \in L^2(\Bbb{R}^d)\::\: (y^2 + (hD_y)^2)^N u \in L^2(\Bbb{R}^d)~~\forall N \in
\Bbb{N}\},$$ the equivalence in Proposition \ref{PropFBIQuadraticEquivalence} reaffirms the well-known fact that $$T(\mathcal{S}(\Bbb{R}^d)) = \{u
\in H_{\Phi}(\Bbb{C}^d;h)\::\: |x|^N u \in L^2_{\Phi}(\Bbb{C}^d;h) ~~\forall N \in \Bbb{N}\}.$$

\end{remark}

\bibliography{MicrolocalBibliography}

\affiliationone{
   Joe Viola\\
   Lund University\\
   Mathematics, Faculty of Science\\
Box 118\\
SE-221 00 LUND\\
Sweden
   \email{jviola@maths.lth.se}}

\end{document}